\documentclass[12pt]{article}

\usepackage{geometry}
\usepackage[left]{lineno}

\usepackage{amsmath, amssymb, amsthm, amsfonts}
\usepackage{mathrsfs, latexsym, euscript, bm, cancel, braket}

\usepackage[all]{xy}
\usepackage{tikz-cd}
\usepackage{graphicx}

\usepackage[inline]{enumitem}
\usepackage{framed, multicol}
\usepackage{comment}
\usepackage{program}
\catcode`\|=12\relax

\usepackage[
 bookmarks=true,
 bookmarksdepth=tocdepth,
 bookmarksnumbered=true,
 colorlinks=false,
 pdftitle={},
 pdfsubject={},
 pdfauthor={},
 pdfkeywords={}
]{hyperref}

\theoremstyle{plain}
\newtheorem{theorem}{\underline{\bf Theorem}}[section]
\newtheorem{lem}[theorem]{\underline{\bf Lemma}}
\newtheorem{cor}[theorem]{\underline{\bf Corollary}}

\newtheorem{prop}[theorem]{\underline{\bf Proposition}}

\theoremstyle{remark}
\newtheorem{definition}[theorem]{\underline{\bf Definition}}

\newcommand{\DS}{\displaystyle}

\renewcommand{\qed}{\openbox}

\numberwithin{theorem}{section}

\numberwithin{equation}{subsection}

\makeatletter
\renewenvironment{proof}[1][\proofname]{\par
 \pushQED{\qed}%
 \normalfont \topsep6\p@\@plus6\p@\relax
 \trivlist
 \item\relax
 {\itshape
 #1\@addpunct{.}}\hspace\labelsep\ignorespaces
}{%
 \popQED\endtrivlist\@endpefalse
}
\makeatother

\begin{document}

\title{The heat kernel on a complex semisimple Lie group and an integral presentation of the heat kernel on its split real form}
\author{Masafumi Shimada\footnote{Faculty of Mathematics, Kyushu University, 744, Motooka,
 Nishi-ku, Fukuoka-shi, 819-0395, Fukuoka, Japan.\\
 e-mail: \texttt{m.shimada.a90@s.kyushu-u.ac.jp}}}
\date{\empty}
\maketitle

\begin{abstract}
Let $G$ be a connected semisimple Lie group, and $G_0$ be its connected split real form. In this paper, we deduce explicit expressions for the heat kernels $\rho^{G_0}_t$ associated with the Laplace--Beltrami operators $\Delta_{G_0}$ and $\Delta_{G}$ respectively, using the algebra of differential operators on an appropriate homogeneous space. These expressions involve the heat Gaussian and the heat kernel on a maximal compact subgroup. Using these expressions for $\rho^{G_0}_t$ and $\rho^{G}_t$, we derive an integral formula relating the heat kernel $\rho^{G_0}_t$ to $\rho^{G}_t$. In the special case of $G_0=SL(2,\mathbb{R})$, we show that the integral formula of $\rho^{SL(2,\mathbb{R})}$ is expressed in terms of the properties of Tchebycheff polynomials.
\end{abstract}

\tableofcontents

\section{Introduction}

\subsection{The focus of this paper and main problems}
Let $G_0$ be a connected real linear Lie group, and $G$ a complexification of $G$. We denote by $\mathfrak{g}_0$ and $\mathfrak{g}$ the corresponding Lie algebras of $G_0$ and $G$ respectively. We assume that $\mathfrak{g}_0$ admits an inner product $\langle\cdot,\cdot\rangle$. This inner product $\langle\cdot,\cdot\rangle$ on $\mathfrak{g}_0$ then induces a left-$G_0$-invariant Riemannian metric on $G_0$. We extend $\langle\cdot,\cdot\rangle$ on $\mathfrak{g}_0$ to a real-valued inner product on $\mathfrak{g}$ defined by
\begin{equation*}
\langle Z_1+iZ_2,Z'_1+iZ'_2\rangle_{\mathbb{C}}:=\langle Z_1,Z'_1\rangle+\langle Z_2,Z'_2\rangle
\end{equation*}
for $Z_j,Z'_j\in\mathfrak{g}_0$ and $j=1,2$. Then the inner product $\langle\cdot,\cdot\rangle_{\mathbb{C}}$ induces a Riemannian metric on $G$ by left translation. Accordingly, we denote by $\Delta_{G_0}$ and $\Delta_{G}$ the respective Laplace--Beltrami operators on $G_0$ and $G$.

We are now ready to pose the following two fundamental problems on heat kernels. Let $t>0$:
\begin{description}
\item[(A)] Deduce explicit expressions of the heat kernel $\rho^{G_0}_t$ for $\Delta_{G_0}$ on $G_0$ and the heat kernel $\rho^G_t$ for $\Delta_G$.
\item[(B)] Give a relationship between the heat kernels $\rho^{G_0}_t$ on the real group $G_0$ and $\rho_t^G$ on the complex group $G$.
\end{description}
\noindent The main objective of this paper is to deduce an explicit formula for $\rho^{G_0}_t$ and $\rho^G_t$ in Problem (A) if $\mathfrak{g}_0$ is a split real form of $\mathfrak{g}$, by using spherical functions on $G_0$, those on $G$ and the heat kernel on a maximal compact subgroup.

\subsection{Background}
We explain our motivation for posing Problems (A) and (B) in terms of heat kernel transforms and a reduction of spherical functions in the real case to those in the complex case.\\

Problems (A) and (B) play an important role in the construction of the Segal--Bargmann--Hall transform for $G_0$ with respect to the heat kernel measures if $G_0$ is a Euclidean space or a connected compact Lie group. For $G_0=\mathbb{R}^N$, the classical Segal--Bargmann transform $B_t$ for $t>0$ is a Hilbert space isometry between $L^2(\mathbb{R}^N)$ and the Segal--Bargmann space of square integrable complex analytic functions on $G=\mathbb{C}^N$ with respect to the Gaussian measure $\rho^{\mathbb{C}^N}_t(z)dz$. The Segal--Bargmann transform $B_t$ sends an initial condition $f\in C^0(\mathbb{R}^N)$ to the analytic continuation of $e^{t\Delta_{\mathbb{R}^N}}f$, and intertwines the canonical commutation relations in those two Hilbert spaces. It is worth pointing out that the Segal--Bargmann transform for $\mathbb{R}^N$ appears in classical results of Euclidean harmonic analysis, e.g., analytic continuation of suitable functions $f$ on $G_0=\mathbb{R}^N$ to appropriate domains in $G=\mathbb{C}^N$, and characterization of the holomorphic extension in terms of $f$ itself or its transforms.

The Segal--Bargmann transform for $\mathbb{R}^N$ has been generalized to many other groups in the spirit of constructions of geometric quantization \cite[Subsection 3.2]{Hall01}. In the case when $G_0$ is a connected Lie group of compact type, i.e., it is isomorphic to a product of $\mathbb{R}^N$ and a connected compact Lie group, complete results \cite{Hall94,Hall01} have been obtained for the corresponding transform called the Segal--Bargmann--Hall transform. If $G_0$ is a connected real linear Lie group and not of compact type, there are sporadic earlier results on the generalized Segal--Bargmann transform for $G_0$. For instance, Kr\"otz--Thangavelu--Xu \cite{KTX05,KTX08} studies the image of the heat kernel transform for a Heisenberg group, and Segal--Bargmann transforms for several motion groups are given in \cite{Sen16} and its references. Among previous results for such noncompact groups $G_0$, one nice feature is that we obtain the special formulae for the heat kernel on $G_0$ and that of $G$, which reflect a structure of the heat kernel governed by special functions. However these special formulae are not available for general pairs of noncompact groups. Thus, for noncompact connected real linear Lie groups $G_0$, calculating the heat kernels $\rho^{G_0}_t$ and $\rho^{G}_t$ explicitly in Problem (A), and investigating the relationship between the two heat kernel measures in Problem (B) can be an initial attempt to explain the existence of a unitary intertwining map in a suitable situation.\\

Problems (A) and (B) is to find an analogous reduction of harmonic analysis on $G_0$ to that on $G$ in relation to Flensted-Jensen's results in spherical harmonic analysis, if $G_0$ is a connected real semisimple Lie group with finite center. In the context of spherical harmonic analysis on a semisimple symmetric space, the Riemannian metric on $G_0$ is replaced by the bi-$G_0$-invariant pseudo-Riemannian metric on $G_0$ induced by the Killing form on $\mathfrak{g}$. Then the heat kernel $\rho^{G_0}_t$ is replaced by the heat Gaussian $g^{G_0}_t$, which is the spherical function defined as the inverse spherical transform of the Gaussian kernel on a Euclidean space. Similarly, we here consider the heat Gaussian $g^{G}_t$ on $G$ instead of the heat kernel $\rho^{G}_t$. Flensted-Jensen establishes a reduction map, denoted by $\mathbf{M}$, from the set of spherical functions on $G_0$ into the set of spherical functions on $G$ that sends the heat Gaussian $g^{G}_t$ to $g^{G_0}_t$. See an explicit formula \eqref{SpecialFunctions:eq4.3.5} for $\mathbf{M}$ in the case of $G_0=SL(2,\mathbb{R})$. We shall emphasize that, for a Riemannian symmetric space $G_0/K_0$ of the noncompact type with a maximal compact subgroup $K_0$ of $G_0$, this map leads to simple direct proofs of results in spherical harmonic analysis such as the spherical Paley--Wiener theorem, the inversion formula and the Plancherel formula. In other words, since spherical functions on $G$ are easy to work with as mentioned in \cite[Section 3]{Fl-J78}, one can reduce the real case to the complex case by this map.

Let us switch to our setting on the Riemannian structures of semisimple $G_0$ and $G$. Mirroring the reduction map $\mathbf{M}$ initiated by Flensted-Jensen \cite{Fl-J78} in spherical harmonic analysis, we expect a correspondence between the set of eigenfunctions of $\Delta_{G_0}$ and the set of eigenfunctions of $\Delta_{G}$ such that the spectral properties of $G_0$ are reduced to the corresponding simpler properties of $G$. However, one should keep in mind that the heat kernels $\rho^{G_0}_t$ and $\rho^{G}_t$ are no longer spherical functions, which precludes the direct application of the reduction map $\mathbf{M}$ to the correspondence of eigenfunctions. To overcome this difficulty, we employ a $K_0$-type decomposition, in which eigenfunctions are expressed as products of radial eigenfunctions and eigenfunctions of the Laplace--Beltrami operator on a maximal compact subgroup $K_0$ of $G_0$. In a similar way, we employ a $K$-type decomposition for $K\subset G$. By initiating Problems (A) and (B), we investigate the possibility of reducing the properties of the heat kernel on $G_0$ to the complex case.

\subsection{Overview of main results}
From now on, we assume that $G$ is a connected complex semisimple Lie group and $G_0$ is its split real form. In Theorem~{\upshape\ref{SpecialFunctions:thm3.0.1}}, we show that the heat kernel $\rho^{G_{\bullet}}_t$ can be expressed by means of the explicit formulae of the heat Gaussian $g^{G_{\bullet}}_t$ and the heat kernel on $K_{\bullet}$. Here $\bullet\in\{0,\emptyset\}$ and we put $G_{\empty}=G$. Let us interpret the statement of Theorem~{\upshape\ref{SpecialFunctions:thm3.0.1}} in terms of the heat kernel transform with respect to the heat kernel measures for $\rho^{G_{\bullet}}_t$. In previous works, to establish the Segal--Bargmann--Hall transform for $K$, one constructs a transform on a dense subspace and extends it to the entire domain. Here the main tool in the construction of the transform is the collection of the matrix entries for finite-dimensional holomorphic representations of its complexification $K_{\mathbb{C}}=G$, which are given by the Peter--Weyl theorem for $K$. Cf. \cite[Section 5.1]{Hall01} for instance. Hence this method for $K$ is independent of any explicit computation of $\rho^{G}_t$ and relies essentially on the compactness of $K$. This reliance may pose an obstacle when computing heat kernel transforms for noncompact groups, such as $G_0$, due to the absence of the corresponding Peter--Weyl type theorem. Note that this obstacle may be related to open problems \cite[Open Problems 1,2]{Hall01} for $G_0$. In this situation, Theorem~{\upshape\ref{SpecialFunctions:thm3.0.1}} suggests that a heat kernel transform for $G_0$ can arise from the replacement of $K$ by the symmetric space $G/K$, via a reduction map \eqref{Chap1:eq1.3.4} sending the heat Gaussian $g^{G}_t$ to $g^{G_0}_t$. Here the map \eqref{Chap1:eq1.3.4} is written explicitly in the equation \eqref{SpecialFunctions:eq4.3.5} if $n=2$. Since the Segal--Bargmann--Hall transform is available for $K$, we hope that the two formulae of $\rho^{G_0}_t$ and $\rho^{G}_t$ lead to a generalized Segal--Bargmann transform for $G_0$ via the reduction map \eqref{Chap1:eq1.3.4}, by switching the roles of $K$ and $G/K$. This provides another approach to the construction of a heat kernel transform for $G_0$ based on an explicit formula for $\rho^{G_{\bullet}}_t$, although we leave the inversion formula for the heat kernel transform to future works.

We move on to the consideration of Problem (B) on a connection between $\rho^{G_0}_t$ and $\rho^{G}_t$. We recall that spherical functions on a complex semisimple Lie group, such as $G$, are in general much simpler objects than their counterparts on a noncompact real form. Theorem~{\upshape\ref{TheHeatKernelOnASemiSimpleGroup:thm4.0.1}} then attempts to construct a counterpart of the reduction map \eqref{SpecialFunctions:eq4.3.5} that reduces harmonic analysis on $G_0$ to that of $G$, although this paper neither specifies any class of functions as its domain nor addresses the construction of an inverse map. Thus a partial answer to Problem (B) is that a desired map sending $\rho^{G}_t$ to $\rho^{G_0}_t$ is given by the combination of the spherical transform \eqref{SpecialFunctions:eq4.3.5} and the integral formula \eqref{TheHeatKernelOnASemiSimpleGroup:eq4.0.1}.

\subsection{Summary of main results}
Let $\bullet\in\{0,\emptyset\}$. We recall that the heat Gaussian $g^{G_{\bullet}}_t$ is defined to be
\begin{equation*}
g^{G_{\bullet}}_t(x)=\frac{1}{|W|}\int_{i\mathfrak{a}^{*}_n}e^{(\langle\lambda,\lambda\rangle-\langle\rho,\rho\rangle)t}\varphi_{\lambda}(x)\frac{d\lambda}{|\mathbf{c}(\lambda)|^2}
\end{equation*}
from the view of spherical harmonic analysis on $G_{\bullet}/K_{\bullet}$. See Proposition~{\upshape\ref{SpecialFunctions:prop2.1.1}} for more details. We introduce a deformation of each spherical function $\varphi_{\lambda}$ with respect to the heat kernel $\rho^{K_{\bullet}}_t$ on $K_{\bullet}$, which is defined by
\begin{equation*}
\int_{K_{\bullet}}e^{(\lambda+\rho)\log\operatorname{Iw}_A(k\sqrt{g(\sigma(g))^{-1}}k^{-1})}\rho^{K_{\bullet}}_t((\operatorname{Iw}_{K_{\bullet}}(k\sqrt{g(\sigma(g))^{-1}}k^{-1})\operatorname{Crt}_{K_{\bullet}}(g)\operatorname{Crt}'_{K_{\bullet}}(g))^{-1})dk
\end{equation*}
for $g\in G_{\bullet}$, in order to attach to right-$K_{\bullet}$-invariant function $\varphi_{\lambda}$ the contribution of the $K_{\bullet}$-cocycle associated with the Iwasawa decomposition of $G_{\bullet}$. Here, $\operatorname{Iw}_{X}(g)$ and $\operatorname{Crt}_{X}(g)$ for $g\in G_{\bullet}$ are determined by the Iwasawa and Cartan decompositions of $G_{\bullet}$ respectively, and $\rho$ denotes the half trace for the corresponding Lie algebra $\mathfrak{g}_{\bullet}$ to $G_{\bullet}$. We now state the main theorem, Theorem~{\upshape\ref{SpecialFunctions:thm3.0.1}}, which presents an explicit formula for the heat kernel $\rho^{G_{\bullet}}_t$.

\theoremstyle{plain}
\newtheorem*{MainResultThm1}{\rm\bf Theorem~{\upshape\ref{SpecialFunctions:thm3.0.1}}}
\begin{MainResultThm1}[Main result]
For $t>0$ and $g\in G_{\bullet}$, the heat kernel $\rho^{G_{\bullet}}_t(g)$ is given by
\begin{equation*}
\begin{split}
\rho^{G_{\bullet}}_t(g)=\frac{1}{|W|}\int_{K_{\bullet}}&\left(\int_{i\mathfrak{a}^{*}_n}e^{(\langle\lambda,\lambda\rangle-\langle\rho,\rho\rangle)\frac{t}{2}}e^{(\lambda+\rho)\log\operatorname{Iw}_A(k\sqrt{g(\sigma(g))^{-1}}k^{-1})}\frac{d\lambda}{|\mathbf{c}(\lambda)|^2}\right)\\
&\times\rho^{K_{\bullet}}_t((\operatorname{Iw}_{K_{\bullet}}(k\sqrt{g(\sigma(g))^{-1}}k^{-1})\operatorname{Crt}_{K_{\bullet}}(g)\operatorname{Crt}'_{K_{\bullet}}(g))^{-1})dk.
\end{split}
\end{equation*}
\end{MainResultThm1}
\noindent In our proof of Theorem~{\upshape\ref{SpecialFunctions:thm3.0.1}}, the basic tool is the algebra of invariant differential operators on the homogeneous space $G_{\bullet}\times K_{\bullet}/\operatorname{diag}K_{\bullet}$, where $\operatorname{diag}K_{\bullet}$ is the diagonal of $G_{\bullet}\times K_{\bullet}$. The key ingredient is a differential equation on a $G_{\bullet}\times K_{\bullet}$-space of representatives for $G_{\bullet}\times K_{\bullet}/\operatorname{diag}K_{\bullet}$ associated with the Iwasawa decomposition of $G_{\bullet}$. We remark that our approach gives limited information about the joint eigenfunctions of the algebra $\mathbf{D}(G_{\bullet}\times K_{\bullet}/\operatorname{diag}K_{\bullet})$ of invariant differential operators.

By putting together the reduction map \eqref{SpecialFunctions:eq4.3.5}, denoted by $\mathbf{M}$, and an integral expression \eqref{TheHeatKernelOnASemiSimpleGroup:eq4.0.1} of the heat kernel $\rho_t^{K_0}$, Theorem~{\upshape\ref{SpecialFunctions:thm3.0.1}} gives rise to Theorem~{\upshape\ref{TheHeatKernelOnASemiSimpleGroup:thm4.0.1}}, which establishes a certain integral relation between $\rho^{G_0}_t$ and $\rho^G_t$.
\theoremstyle{plain}
\newtheorem*{MainResultThm2}{\rm\bf Theorem~{\upshape\ref{TheHeatKernelOnASemiSimpleGroup:thm4.0.1}}}
\begin{MainResultThm2}
For $t>0$ and $g\in G_0$, the heat kernel $\rho^{G_0}_t(g)$ is expressed as
\begin{equation*}
\begin{split}
\rho_t^{G_0}(g)=c_K\int_{K_0^{\mathbb{C}}}\int_{K_0}\int_{\overline{T^{+}}}\rho_{t/2}^G(ya_gk_gk\tau k^{-1})\delta(\tau)dkd\tau dy+(1-c_K)(\mathbf{M}g_{t/4}^G)(a_g),
\end{split}
\end{equation*}
where $c_K$ is a constant determined by the Haar measure on $K$ and measures associated with the Cartan decomposition of $K$.
\end{MainResultThm2}
\noindent Here, for $g\in G_0$, we put $a_g:=\operatorname{Crt}_{\overline{A^{+}}}(g)$ and $k_g:=\operatorname{Crt}_{K_0}(g)\operatorname{Crt}'_{K_0}(g)$.

\subsection{The organization}
The paper is organized as follows. In Section~{\upshape\ref{TheHeatKernelOnASemiSimpleGroup:sec2}}  we discuss some notation concerning the heat kernel on a connected semisimple Lie group, and recall several fundamental results on the heat Gaussian on the Lie group. Section~{\upshape\ref{TheHeatKernelOnASemiSimpleGroup:sec3}} is devoted to a quadrature of the heat kernel for the Laplace--Beltrami operator on the Lie group. See Theorem~{\upshape\ref{SpecialFunctions:thm3.0.1}}. In Section~{\upshape\ref{TheHeatKernelOnASemiSimpleGroup:sec4}} we relate the heat kernel on a real group to that on its complexification, using Theorem~{\upshape\ref{SpecialFunctions:thm3.0.1}}. In Section~{\upshape\ref{TheHeatKernelOnASemiSimpleGroup:sec5}} we consider the special case of $SL(2,\mathbb{R})$.

\section{Heat kernels and heat Gaussians}\label{TheHeatKernelOnASemiSimpleGroup:sec2}
This section collects some technical results on heat kernels and heat Gaussians, which are needed for the rest of the present paper. We shall summarize briefly the definition and a few basic properties of the heat kernel on a semisimple Lie group. In this section, we provide the preliminary results on the heat kernel and the heat Gaussian required to prove the main theorem.

\subsection{The heat kernel on a Lie group}

Let $G$ be a connected complex semisimple Lie group and $G_0$ be its connected split real form with finite center. Then, for $\bullet\in\{0,\emptyset\}$, the symbol $G_{\bullet}$ is employed when the distinction between real and complex groups is irrelevant, where $G_{\emptyset}=G$. Let $K_{\bullet}$ be a maximal compact subgroup of $G_{\bullet}$. We denote by $e\in G_{\bullet}$ the identity element. Let $\sigma_{\bullet}:\ G_{\bullet}\to G_{\bullet}$ be the Cartan involution corresponding to $K_{\bullet}$. We have the corresponding Lie algebras $\mathfrak{k}_{\bullet}=\operatorname{Lie}(K_{\bullet})$ and $\mathfrak{g}_{\bullet}=\operatorname{Lie}(G_{\bullet})$, where we consider $\mathfrak{g}_{\bullet}$ as a real Lie algebra. Let $\sigma_{\bullet}$ also denote the differential of the Cartan involution on $G_{\bullet}$ such that $\sigma_{\bullet}(\mathfrak{k}_{\bullet})=\mathfrak{k}_{\bullet}$. We denote by $\mathfrak{g}_{\bullet}=\mathfrak{k}_{\bullet}\oplus\mathfrak{p}_{\bullet}$ the Cartan decomposition with respect to $\sigma_{\bullet}$. Let $\exp:\ \mathfrak{g}_{\bullet}\to G_{\bullet}$ be the exponential map of matrices. We denote by $[\cdot,\cdot]$ the Lie bracket of $\mathfrak{g}_{\bullet}$. For $X,Y\in\mathfrak{g}_{\bullet}$, the adjoint action $(\operatorname{ad}X)Y:=[X,Y]$, which is the differential representation of the adjoint representation $\operatorname{Ad}_{G_{\bullet}}(g)Y:=gYg^{-1}$ for $g\in G_{\bullet}$, gives rise to a bilinear form
\begin{equation*}
B(X,Y):=2\operatorname{tr}_{\mathbb{R}}
(\operatorname{ad}X\circ\operatorname{ad}Y).
\end{equation*}
By translation, the inner product $(X,Y)\mapsto-B(X,\sigma_{\bullet}(Y))$ defines a left-$G_{\bullet}$-invariant and right-$K_{\bullet}$-invariant Riemannian metric $(g_{ij})$ on $G_{\bullet}$. A Haar measure $dg$ on $G_{\bullet}$ is defined as the volume form of $(g_{ij})$. Let $\omega_{G_{\bullet}}$ be the Casimir operator and $\Delta_{G_{\bullet}}$ be the Laplace--Beltrami operator, which are induced by $B$ and $(g_{ij})$ respectively. Thus, the heat operator is $L_{\Delta}=\partial_t-\Delta_{G_{\bullet}}$.

Let $B_{\mathfrak{k}_{\bullet}}=2^{-1}\left.B\right|_{\mathfrak{k}_{\bullet}}$ be the Killing form on $\mathfrak{k}_{\bullet}$, and let $\omega_{K_{\bullet}}$ be the Casimir operator associated with $B_{\mathfrak{k}_{\bullet}}$. Then $B_{\mathfrak{k}_{\bullet}}$ induces a Riemannian structure on $K_{\bullet}$ by a similar argument. Hence let $\Delta_{K_{\bullet}}$ denote the Laplace--Beltrami operator on $K_{\bullet}$ induced by $B_{\mathfrak{k}_{\bullet}}$. We have $\Delta_{K_{\bullet}}=-\omega_{K_{\bullet}}$ and $\Delta_{G_{\bullet}}=\omega_{G_{\bullet}}-\omega_{K_{\bullet}}$.\\

We discuss the existence and uniqueness of the heat kernel on $G_{\bullet}$. See \cite{Shima26Thesis}[Theorem 1.1.1] in the case of special linear groups. It is shown in \cite{Dodz83} that the fundamental solution of the heat equation for $L_{\Delta}$ exists uniquely on the complete Riemannian manifold $G_{\bullet}$ with Ricci curvature bounded from below. Recall that, by the definition of the left-$G_{\bullet}$-invariant metric $(g_{ij})$, $G_{\bullet}$ acts transitively on itself by isometries. Hence, using \cite[Chapter IV, Theorem 4.5]{KN63} on the geodesical completeness of a homogeneous Riemannian manifold, and applying the Hopf--Rinow theorem, we find that $G_{\bullet}$ is a complete Riemannian manifold. Moreover, putting together the Cartan decomposition $\mathfrak{g}_{\bullet}=\mathfrak{k}_{\bullet}\oplus\mathfrak{p}_{\bullet}$, we deduce from the homogeneity of $G_{\bullet}$ that the Ricci curvature of $G_{\bullet}$ is bounded from below. See \cite{Shima26Thesis}[Chapter 1, Section 1.1, pp. 31--33] for an estimate of the Ricci curvature in the case of $SL(n,\mathbb{K})$ for $\mathbb{K}=\mathbb{R},\mathbb{C}$. Thus we obtain Theorem~{\upshape\ref{SpecialFunctions:thm2.1.1}} on the existence and the uniqueness of the fundamental solution for $L_{\Delta}$.
\begin{theorem}[{Cf. \cite{Dodz83}}]\label{SpecialFunctions:thm2.1.1}
There exists a unique fundamental solution associated with $\Delta_{G_{\bullet}}$ on $G_{\bullet}$.
\end{theorem}
\noindent Therefore, we define the heat kernel on $G_{\bullet}$ for $L_{\Delta}$ in Definition~{\upshape\ref{SpecialFunctions:def2.1.1}}. Since the fundamental solution is left-$G_{\bullet}$ and right-$K_{\bullet}$-invariant, we get the heat kernel at the identity $\rho^{G_{\bullet}}_t$ of a single variable. Cf. \cite[Theorem 9.12]{Gri09}. In particular, since the heat kernel on $G_{\bullet}$ is symmetric with respect to its spatial variables, we have $\rho^{G_{\bullet}}_t(g^{-1})=\rho^{G_{\bullet}}_t(g)$ for $g\in G_{\bullet}$.

\begin{definition}\label{SpecialFunctions:def2.1.1}
The heat kernel $\rho^{G_{\bullet}}_t(g)$ of the heat equation for $\Delta_{G_{\bullet}}$ is the fundamental solution of the problem \eqref{SpecialFunctions:eq2.1.1}:
\begin{equation}\label{SpecialFunctions:eq2.1.1}
L_{\Delta}\rho^{G_{\bullet}}_t(g)=0,\ \|\rho^{G_{\bullet}}_t*f-f\|_{L^2}\to0
\end{equation}
as $t\to0+$ for compactly supported functions $f\in C^{\infty}_c(G_{\bullet})$, where $\|f\|_{L^2}$ is the $L^2$-norm of $f$ in $L^2(G_{\bullet},dg)$ and $\rho^{G_{\bullet}}_t*f(x):=\int_{G_{\bullet}}\rho^{G_{\bullet}}_t(g)f(xg^{-1})dg$ denotes the convolution product.
\end{definition}
\noindent Note that $\rho^{G_{\bullet}}_t(g)$ is $C^{\infty}$-smooth with respect to both $t>0$ and $g\in G_{\bullet}$, and it also belongs to $L^2(G_{\bullet},dg)$ for $t>0$. Cf. \cite[Theorem 7.20]{Gri09}.\\

Recall that there uniquely exists a fundamental solution for $\Delta_{K_{\bullet}}$ on the compact group $K_{\bullet}$. In a manner similar to Definition~{\upshape\ref{SpecialFunctions:def2.1.1}}, let $\rho^{K_{\bullet}}_t=\rho^{K_{\bullet}}_t(k)$ denote the heat kernel on $K_{\bullet}$. Theorem~{\upshape\ref{SpecialFunctions:thm2.1.2}} states that $\rho^{K_{\bullet}}_t$ has a series expansion with respect to the characters of the irreducible unitary representations of $K_{\bullet}$. Take an irreducible unitary representation $(\tau,V_{\tau})$ of $K_{\bullet}$. We denote by $d\tau$ the infinitesimal action of $\mathfrak{k}_{\bullet}$. Let $\chi_{\tau}(k):=\operatorname{tr}_{V_{\tau}}\tau(k)$ be the character of $\tau$ for $k\in K_{\bullet}$. We denote by $\hat{K}_{\bullet}$ the set of isomorphism classes of irreducible unitary representations of $K_{\bullet}$.

\begin{theorem}[cf. {\cite[Chapter 2, Theorem 1]{Stein70}}]\label{SpecialFunctions:thm2.1.2}
The heat kernel $\rho^{K_{\bullet}}_t$ exists uniquely as the fundamental solution to the heat equation for the heat operator $\partial_t-\Delta_{K_{\bullet}}$, and is given by
\begin{equation}\label{SpecialFunctions:eq2.1.5}
\rho^{K_{\bullet}}_t(k)=\sum_{[(\tau,V_{\tau})]\in\hat{K}_{\bullet}}(\operatorname{dim}V_{\tau})e^{-\lambda_{\tau}t}\chi_{\tau}(k).
\end{equation}
\end{theorem}
\noindent Here, $\lambda_{\tau}$ is a non-negative real number determined by Schur's lemma since $\tau$ is irreducible, which satisfies the equation
\begin{equation*}
\sum_m d\tau(Z_m)\circ d\tau(Z_m)=-\lambda_{\tau}\operatorname{id}_{V_{\tau}}
\end{equation*}
for any orthonormal basis $\{Z_m\}$ of $\mathfrak{k}_{\bullet}$ with respect to $-B_{\mathfrak{k}_{\bullet}}$. Note that $\lambda_{\tau}$ is independent of the choice of the basis $\{Z_m\}$.

\subsection{The heat Gaussian on a Riemnnain symmetric space}
Let $L_{\omega}=\partial_t-2\omega_{G_{\bullet}}$. It follows from Proposition~{\upshape\ref{SpecialFunctions:prop2.1.1}} and Lemma~{\upshape\ref{SpecialFunctions:lem2.1.1}} that a solution, denoted by $g^{G_{\bullet}}_t$, of the problem analogous to \eqref{SpecialFunctions:eq2.1.1} for the operator $L_{\omega}$ is determined as the inverse spherical transform of the Gaussian on a Euclidean space. Hence $g^{G_{\bullet}}_t$ is called the heat Gaussian. Cf. \cite[Section 3]{Gan68} and \cite[Chapter 2, Section 2]{JL09}.

Before stating the expression \eqref{SpecialFunctions:eq2.1.2} of the heat Gaussian $g^{G_{\bullet}}_t$ in Proposition~{\upshape\ref{SpecialFunctions:prop2.1.1}}, we establish the necessary notation. Recall the Cartan decomposition $\mathfrak{g}_{\bullet}=\mathfrak{k}_{\bullet}\oplus\mathfrak{p}_{\bullet}$ with respect to $\sigma$. Let $\mathfrak{a}_{\bullet}\subset\mathfrak{p}_{\bullet}$ be a maximal abelian subspace. Note that $\mathfrak{a}_0=\mathfrak{a}$. Let $\mathfrak{a}^{*}_{\bullet}$ be the dual of $\mathfrak{a}_{\bullet}$, $\Sigma_{\bullet}\subset\mathfrak{a}^{*}_{\bullet}$ be the corresponding set of restricted roots and $W=W(\mathfrak{g}_{\bullet},\mathfrak{a}_{\bullet})$ be the Weyl group. Here, $|W|$ denotes the order of $W$. Let $\langle\cdot,\cdot\rangle$ be a bilinear form on $\mathfrak{a}^{*}_{\bullet}\oplus i\mathfrak{a}^{*}_{\bullet}$ induced by the Killing form $B$. Fix a Weyl chamber $\mathfrak{a}^{+}_{\bullet}\subset\mathfrak{a}_{\bullet}$; we denote by $\Sigma^{+}_{\bullet}$ the corresponding set of positive restricted roots, and by $\mathfrak{n}:=\sum_{\alpha\in\Sigma^{+}_{\bullet}}(\mathfrak{g}_{\bullet})_{\alpha}$ the nilpotent algebra spanned by the root subspaces $(\mathfrak{g}_{\bullet})_{\alpha}:=\{X\in\mathfrak{g}_{\bullet}\ ;\ (\operatorname{ad}H)X=\alpha(H)X\ (H\in\mathfrak{a}_{\bullet})\}$. Let $m_{\alpha}:=\dim(\mathfrak{g}_{\bullet})_{\alpha}$ and let $\rho:=2^{-1}\sum_{\alpha\in\Sigma^{+}_{\bullet}}m_{\alpha}\alpha$ be the half trace. If $A_{\bullet}$ and $N$ denote the respective analytic subgroups of $G_{\bullet}$ with Lie algebras $\mathfrak{a}_{\bullet}$ and $\mathfrak{n}$, we have the Iwasawa decomposition $G_{\bullet}=NA_{\bullet}K_{\bullet}$. Let $\overline{A^{+}_{\bullet}}\subset A_{\bullet}$ be the subset corresponding to the closure of $\mathfrak{a}^{+}_{\bullet}$, then we have the Cartan decomposition $G_{\bullet}=K_{\bullet}\overline{A^{+}_{\bullet}}K_{\bullet}$. Let $\operatorname{Iw}_A,\operatorname{Iw}_K$ be the respective projections of $G_{\bullet}=NA_{\bullet}K_{\bullet}$ on $A_{\bullet}$ and $K_{\bullet}$. For each $g\in G_{\bullet}$, it is known that, given a Cartan decomposition $g=k_1ak_2\in K_{\bullet}\overline{A^{+}_{\bullet}}K_{\bullet}$, the element $a\in\overline{A^{+}}$ is uniquely determined. Hence we denote by $\operatorname{Crt}_{\overline{A^{+}}}$ the projection of $G_{\bullet}=K_{\bullet}\overline{A^{+}_{\bullet}}K_{\bullet}$ on $\overline{A^{+}_{\bullet}}$. For $a\in A_{\bullet}$, we denote by $\log a$ the unique element $H\in\mathfrak{a}_{\bullet}$ such that $a=\exp H$. For each $g\in G_{\bullet}$ with $\operatorname{Crt}_{\overline{A^{+}}}(g)\neq e$, given a decomposition $g=k_1\operatorname{Crt}_{\overline{A^{+}}}(g)k_2\in K_{\bullet}\overline{A^{+}_{\bullet}}K_{\bullet}$, we put $\operatorname{Crt}_K(g):=k_1$ and $\operatorname{Crt}'_K(g):=k_2$. If $k\in K_{\bullet}$, we put $\operatorname{Crt}_K(k):=k$ and $\operatorname{Crt}'_K(k):=e$. Then we shall verify that, if $k_1ak_2=k'_1ak'_2\in K_{\bullet}\overline{A^{+}_{\bullet}}K_{\bullet}$, there exists $x\in Z_{K_{\bullet}}(a)$ such that $k'_1=k_1x$ and $k'_2=x^{-1}k_2$, where $Z_{K_{\bullet}}(a)$ denotes the centralizer of $a\in\overline{A^{+}_{\bullet}}$ in $K_{\bullet}$. Without loss of generality, we can assume that $k'_1=k'_2=e$, or equivalently, $k_1ak_2=a$. Take the unique $H\in\mathfrak{a}_{\bullet}$ such that $a=\exp H$. It follows from $\sigma(\exp (-H))=\exp H$ that
\begin{equation*}
\begin{split}
&\exp(2H)=(\exp H)\sigma(\exp (-H))=k_1(\exp H)k_2\sigma(k^{-1}_2\exp (-H)k^{-1}_1)\\
&=k_1(\exp (2H))k^{-1}_1=\exp(2k_1Hk^{-1}_1),
\end{split}
\end{equation*}
which yields $k_1Hk^{-1}_1=H$. Hence we get $k_1ak^{-1}_1=\exp(k_1Hk^{-1}_1)=\exp H=a$, that is, $k_1\in Z_{K_{\bullet}}(a)$. Thus we deduce from $a=k_1ak_2=a(k_1k_2)$ that $k_2=k^{-1}_1\in Z_{K_{\bullet}}(a)$. Therefore we obtain the fact that, given $g\in G_{\bullet}$, the expressions $\operatorname{Crt}_K(g)\operatorname{Crt}'_K(g)$ and $\operatorname{Crt}_K(g)\operatorname{Crt}_{\overline{A^{+}}}(g)(\operatorname{Crt}_K(g))^{-1}$ are well-defined.

We denote by
\begin{equation*}
\varphi_{\lambda}(g)=\int_{K_{\bullet}}e^{(\rho+\lambda)(\log\operatorname{Iw}_A(kg))}dk
\end{equation*}
the spherical function corresponding to $\lambda\in\mathfrak{a}^{*}_{\bullet}\oplus i\mathfrak{a}^{*}_{\bullet}$. Recall that the family $\{\varphi_{\lambda}\}$, parametrized by $\lambda\in\mathfrak{a}^{*}_{\bullet}\oplus i\mathfrak{a}^{*}_{\bullet}$, forms a system of joint eigenfunctions of the algebra $\mathbf{D}(G_{\bullet}/K_{\bullet})$ of differential operators on the symmetric space $G_{\bullet}/K_{\bullet}$ that are invariant under the left translations from $G_{\bullet}$, and that we have
\begin{equation*}
\omega_{G_{\bullet}}\varphi_{\lambda}=(\langle\lambda,\lambda\rangle-\langle\rho,\rho\rangle)\varphi_{\lambda}.
\end{equation*}
Let $\mathbf{c}=\mathbf{c}(\lambda)$ denote Harish-Chandra's $\mathbf{c}$-function on $i\mathfrak{a}^{*}_{\bullet}$. It is known that $\mathbf{c}(\lambda)$ does not vanish on $i\mathfrak{a}^{*}_{\bullet}$ and satisfies a product formula written as
\begin{equation*}
\mathbf{c}(\lambda)=C_0\prod_{\alpha}\frac{2^{-\langle\lambda,\alpha_0\rangle}{\boldsymbol \Gamma}(\langle\lambda,\alpha_0\rangle)}{{\boldsymbol \Gamma}\left(\frac{m_{\alpha}}{4}+\frac{1}{2}+\frac{\langle\lambda,\alpha_0\rangle}{2}\right){\boldsymbol \Gamma}\left(\frac{m_{\alpha}}{4}+\frac{m_{2\alpha}}{2}+\frac{\langle\lambda,\alpha_0\rangle}{2}\right)},
\end{equation*}
where ${\boldsymbol \Gamma}$ is the Gamma function, $\alpha$ runs over $\Sigma^{+}_{\bullet}$, $\alpha_0:=\alpha/\langle\alpha,\alpha\rangle$, and the constant $C_0$ is given by $\mathbf{c}(-i\rho)=1$.

\begin{prop}[cf. {\cite[Proposition 3.1]{Gan68}}]\label{SpecialFunctions:prop2.1.1}
For $x\in G_{\bullet}$, the heat Gaussian
\begin{equation}\label{SpecialFunctions:eq2.1.2}
g^{G_{\bullet}}_t(x)=\frac{1}{|W|}\int_{i\mathfrak{a}^{*}_{\bullet}}e^{(\langle\lambda,\lambda\rangle-\langle\rho,\rho\rangle)t}\varphi_{\lambda}(x)\frac{d\lambda}{|\mathbf{c}(\lambda)|^2}
\end{equation}
is a continuous function that satisfies $L_{\omega}g^{G_{\bullet}}_t(x)=0$, and $\|g^{G_{\bullet}}_t*f-f\|_{L^2}\to0$ as $t\to0+$ for spherical $f\in C^{\infty}_c(G_{\bullet})$.
\end{prop}
\noindent We recall several properties of the heat Gaussian $g^{G_{\bullet}}_t$. Given $x\in G_{\bullet}$, we use the notation
\begin{equation*}
\|x\|:=\sqrt{\langle\log\operatorname{Crt}_{\overline{A^{+}}}(x),\log\operatorname{Crt}_{\overline{A^{+}}}(x)\rangle},\end{equation*}
and call $\|x\|$ the polar height of $x$. If $\operatorname{dist}=\operatorname{dist}(xK_{\bullet},yK_{\bullet})$ denotes the distance on $G_{\bullet}/K_{\bullet}$ associated with the Cartan involution $\sigma$, one finds that
\begin{equation*}
\operatorname{dist}(xK_{\bullet},yK_{\bullet})=\|x^{-1}y\|.
\end{equation*}

\begin{lem}[cf. {\cite[Chapter X, Section 7, and Chapter XII, Corollary 5.2]{JL01}}]\label{SpecialFunctions:lem2.1.1}

\begin{enumerate}
\item For each $t>0$, $g^{G_{\bullet}}_t$ is positive and satisfies
\begin{equation}\label{SpecialFunctions:eq2.1.3}
\int_{G_{\bullet}}g^{G_{\bullet}}_t(x)dx=1.
\end{equation}
\item Given $\delta>0$, then we have
\begin{equation}\label{SpecialFunctions:eq2.1.4}
\lim_{t\to0+}\int_{\|x\|\geq\delta}g^{G_{\bullet}}_t(x)dx=0.
\end{equation}
\end{enumerate}
\end{lem}
The heat Gaussian $g^{G_{\bullet}}_t$ induces the fundamental solution for the Laplace--Beltrami operator on the symmetric space $G_{\bullet}/K_{\bullet}$. More precisely, since $g^{G_{\bullet}}_t$ is bi-$K_{\bullet}$-invariant, we obtain the function on $(0,\infty)\times G_{\bullet}/K_{\bullet}\times G_{\bullet}/K_{\bullet}$
\begin{equation*}
\mathbf{K}^{G_{\bullet}/K_{\bullet}}(t,xK_{\bullet},yK_{\bullet}):= g^{G_{\bullet}}_t(x^{-1}y).
\end{equation*}
It is verified in \cite{Shima26Thesis}[Chapter 1, Section 1.1, pp. 36--39] that $\mathbf{K}^{G_{\bullet}/K_{\bullet}}$ is the heat kernel for the Laplace--Beltrami operator on $G_{\bullet}/K_{\bullet}$, where one key ingredient is Lemma~{\upshape\ref{Chap1:lem1.1.2}} for the symmetric space $G_{\bullet}/K_{\bullet}$.

\begin{lem}[{Cf. \cite[Lemma 9.2]{Gri09}}]\label{Chap1:lem1.1.2}
Let $M$ be a Riemannian manifold. We denote by $\mathbf{K}^M$ the heat kernel associated with the Laplace--Beltrami operator on $M$. Let $u(t,x)$ be a smooth non-negative function on $(0,\infty)\times M$ satisfying
\begin{equation*}
\int_M u(t,x)d\mu(x)\leq1
\end{equation*}
for $t>0$, where $d\mu$ denotes the volume form on $M$. Then the following conditions (i)--(iii) are equivalent:
\begin{enumerate}
\item $u(t,\bullet)\to\delta_y$ as $t\to0+$ in the distribution topology on $M$, where $\delta_y$ denotes the Dirac delta function at $y\in M$.
\item For any open set $U$ containing $y$,
\begin{equation*}
\lim_{t\to0+}\int_U u(t,x)d\mu(x)=1.
\end{equation*}
\item For any bounded continuous function $f$ on $M$,
\begin{equation*}
\int_M u(t,x)d\mu(x)\to f(y)
\end{equation*}
as $t\to0+$.
\end{enumerate}
\end{lem}

We recall that Flensted-Jensen \cite{Fl-J78} establishes a reduction of spherical functions on $G_0$ to elementary spherical functions on $G$, which relates the heat Gaussian $g_t^{G_0}$ to $g_t^G$. We proceed as in \cite{Fl-J78}. Since the maximal abelian subspace $\mathfrak{a}$ is the same as $\mathfrak{a}_0$ in this setting, we omit the subscript on $\mathfrak{a}_0$, and hence also on $(\mathfrak{a}_0)^{*}$ and $A_0$. We set $\mathfrak{a}^{*}_{\mathbb{C}}:=\mathfrak{a}^{*}\oplus i\mathfrak{a}^{*}$. Let $\varphi_{\lambda}$ denote the bi-$K_0$-invariant functions on $G_0$ for $\lambda\in\mathfrak{a}^{*}_{\mathbb{C}}$, and $\Phi_{\Lambda}$ denote the bi-$K$-invariant functions on $G$ for $\Lambda\in\mathfrak{a}^{*}_{\mathbb{C}}$.

Let us use the same symbol $\sigma$ for the Cartan involution on $G$ and its corresponding involution on $\mathfrak{g}$, because they are fundamentally the same geometric transformation viewed at different scales. Similarly, we use the same notation $\sigma_0$ for the Cartan involution on $\mathfrak{g}_0$ corresponding to the Cartan involution on $G_0$. We denote by $\tilde{\sigma}_0$ the complex linear extension of $\sigma_0$ to $\mathfrak{g}$, and let $\mathfrak{k}^{\mathbb{C}}_0$ be the fixed point set of $\tilde{\sigma}_0$. Then, as a generalization of the Cartan decomposition $G=K\overline{A^{+}}K=G^{\sigma}\overline{A^{+}}G^{\sigma}$, we state the generalized Cartan decomposition \eqref{Chap1:eq1.3.1} of $G$.
\begin{theorem}[{Cf. \cite[Theorem 4.1]{Fl-J78} and \cite[Section 3, Theorem 10]{Ross79}}]\label{Chap1:thm1.3.4}
The group $G$ admits a decomposition
\begin{equation}\label{Chap1:eq1.3.1}
G=G^{\tilde{\sigma}_0}\overline{A^{+}}G^{\sigma},
\end{equation}
and, for each $g\in G$, there exists a unique $a\in\overline{A^{+}}$ such that $g\in G^{\tilde{\sigma}_0} a G^{\sigma}$.
\end{theorem}
\noindent We remark that, in the geometric context, Rossmann \cite{Ross79} proves the generalized Cartan decomposition of a semisimple symmetric space into orbits of a maximal compact subgroup of the group of displacements.

Via the generalized Cartan decomposition \eqref{Chap1:eq1.3.1}, the analysis of the double coset space $K_0\times K_0\backslash G_0\times G_0/\operatorname{diag}G_0$ amounts to Theorem~{\upshape\ref{Chap1:thm1.3.5}}.

\begin{theorem}[{Cf. \cite[Theorem 5.5]{Fl-J78}}]\label{Chap1:thm1.3.5}
If $\lambda\in \mathfrak{a}^{*}_{\mathbb{C}}$, define $\Lambda_{\lambda}\in \mathfrak{a}^{*}_{\mathbb{C}}$ by
\begin{equation*}
\Lambda_{\lambda}+i\rho=2(\lambda+i\rho_0).
\end{equation*}
Then, for $x\in G$, we have
\begin{equation}\label{Chap1:eq1.3.2}
\Phi_{\Lambda_{\lambda}}(x)=\int_K\varphi_{\lambda}(kx(\sigma(kx))^{-1})dk.
\end{equation}
\end{theorem}
\noindent We denote by $\pi_0$ the product of positive restricted roots in $\Sigma^{+}_0$ for the pair $(\mathfrak{g}_0,\mathfrak{a}_0)$. If the Haar measure $dh$ on $K^{\mathbb{C}}_0$ is suitably normalized, the inverse relation associated with the formula \eqref{Chap1:eq1.3.2} is given by
\begin{equation}\label{Chap1:eq1.3.3}
\varphi_{\lambda}(x(\sigma(x))^{-1})=|\mathbf{c}(\lambda)|^2|\pi_0(\lambda)|^2\int_{K^{\mathbb{C}}_0}\Phi_{2\lambda}(hx)dh
\end{equation}
for $\lambda\in i\mathfrak{a}^{*}$ and $x\in G_0$. Cf. \cite[Equation (1.4) and Corollary 7.4]{Fl-J78}. Hence, in particular, the Plancherel measure for $G_0/K_0$ is presented by
\begin{equation*}
|\mathbf{c}(\lambda)|^{-2}=|\pi_0(\lambda)|^2\int_{K^{\mathbb{C}}_0}\Phi_{2\lambda}(h)dh
\end{equation*}
for $\lambda\in i\mathfrak{a}^{*}$ in terms of spherical functions on $G$. Here $\mathbf{c}$ denotes Harish-Chandra's $\mathbf{c}$-function for $G_0$. Cf. \cite[Equations (1.3) and (7.6)]{Fl-J78}. It is known that the formulae \eqref{Chap1:eq1.3.2} and \eqref{Chap1:eq1.3.3} describing the correspondence between spherical functions on $G$ and those on $G_0$ induce a convolution-preserving bijection 
\begin{equation}\label{Chap1:eq1.3.4}
\mathbf{M}:\ \mathscr{S}^2(K\backslash G/K)\to \mathscr{S}^2(K_0\backslash G_0/K_0),
\end{equation}
\noindent which is described explicitly as the map \eqref{SpecialFunctions:eq4.3.5} for $G_0=SL(2,\mathbb{R})$. Here $\mathscr{S}^2(K_0\backslash G_0/K_0)\subset C^{\infty}(K_0\backslash G_0/K_0)$ denotes the $L^2$-spherical-Schwartz space on $G_0$ defined by
\begin{equation*}
\begin{split}
&\mathscr{S}^2(K_0\backslash G_0/K_0)\\
&:=\{f\ ;\ 1\leq {}^{\forall}N\in\mathbb{Z},\ {}^{\forall}D\in\mathbf{D}(G_0),\ \sup_{x\in G_0}|Df(x)|(1+\|x\|)^N(\left.\varphi_{\lambda}\right|_{\lambda=0}(x))^{-1}<\infty\}
\end{split}
\end{equation*}
with the algebra $\mathbf{D}(G_0)$ of left-$G_0$-invariant differential operators on $G_0$. Cf. \cite[Theorem 7.5]{Fl-J78}, and cf. \cite[Section 7, Remark, p. 135]{Fl-J78} for the equivalence between the various definitions of the $L^2$-spherical-Schwartz spaces. In particular, we have
\begin{equation}\label{Chap1:eq1.3.7}
g^{G_0}_{t/2}=\mathbf{M}g^G_{t/4}\in\mathscr{S}^2(K_0\backslash G_0/K_0)
\end{equation}
for $t>0$. Cf. \cite[Equation (6.20)]{Fl-J78}.

\section{An explicit formula for the heat kernel on a semisimple Lie group}\label{TheHeatKernelOnASemiSimpleGroup:sec3}
In this section, we prove Theorem~{\upshape\ref{SpecialFunctions:thm3.0.1}} concerning a quadrature of the heat kernel $\rho^{G_{\bullet}}_t$ for the Laplace--Beltrami operator $\Delta_{G_{\bullet}}$ on $G_{\bullet}$, where $\bullet\in\{0,\emptyset\}$.
\begin{theorem}\label{SpecialFunctions:thm3.0.1}
Let $\bullet\in\{0,\emptyset\}$ and $g\in G_{\bullet}$. We define
\begin{equation*}
\begin{split}
&\sqrt{g(\sigma_{\bullet}(g))^{-1}}:=\operatorname{Crt}_{K_{\bullet}}(g)\operatorname{Crt}_{\overline{A^{+}}}(g)(\operatorname{Crt}_{K_{\bullet}}(g))^{-1}.
\end{split}
\end{equation*}
Then for $t>0$, the heat kernel $\rho^{G_{\bullet}}_t(g)$ is given by
\begin{equation}\label{SpecialFunctions:eq3.0.1}
\begin{split}
\rho^{G_{\bullet}}_t(g)=\frac{1}{|W|}\int_{K_{\bullet}}&\left(\int_{i\mathfrak{a}^{*}_n}e^{(\langle\lambda,\lambda\rangle-\langle\rho,\rho\rangle)\frac{t}{2}}e^{(\lambda+\rho)\log\operatorname{Iw}_A(k\sqrt{g(\sigma_{\bullet}(g))^{-1}}k^{-1})}\frac{d\lambda}{|\mathbf{c}(\lambda)|^2}\right)\\
&\times\rho^{K_{\bullet}}_t((\operatorname{Iw}_{K_{\bullet}}(k\sqrt{g(\sigma_{\bullet}(g))^{-1}}k^{-1})\operatorname{Crt}_{K_{\bullet}}(g)\operatorname{Crt}'_{K_{\bullet}}(g))^{-1})dk.
\end{split}
\end{equation}
\end{theorem}
\noindent In the following subsections, we only deduce an expression of the heat kernel $\rho_t^G=\rho_t^{G_{\emptyset}}$, since a similar argument works for the heat kernel $\rho_t^{G_0}$ by adding the subscript zero to the objects associated with $G$. More precisely, we do not use any property specific to a complex group in this section.

\subsection{Invariant differential operators on a homogeneous space}
Recall that $\mathbf{D}(G\times K)$ 
is linearly isomorphic to the symmetric algebra on the vector space $\mathfrak{g}\oplus\mathfrak{k}$. Let $H$ be the image of the diagonal embedding $K\hookrightarrow G\times K$, and $\mathfrak{h}\subset\mathfrak{g}\oplus\mathfrak{k}$ be the Lie algebra of $H$. Let $G\times K\curvearrowleft K$ with $(g,k)\cdot k':=(gk',kk')$ the diagonal action, and let $G\times K/H\to G/K$ be a homogeneous principal $K$-bundle. We denote by $\pi:\ G\times K\to G\times K/H$ the natural projection.

We recall the algebra of the left-$G\times K$-invariant differential operators on the homogeneous space $G\times K/H$, denoted by $\mathbf{D}(G\times K/H)$. Let us take a basis $\{X_{K,i}\}\subset\mathfrak{k}$ such that $\omega_K=-2\sum_i\tilde{X}^2_{K,i}$ on $K$, where $\tilde{X}_{K,i}$ is the left invariant vector field on $K$ associated with $X_{K,i}$. Via the inclusion $\mathfrak{k}\hookrightarrow\mathfrak{g}$ induced by $K\hookrightarrow G$, we denote by $X_i$ the image of $X_{K,i}$. Under the Cartan decomposition $\mathfrak{g}=\mathfrak{k}\oplus\mathfrak{p}$, we take a basis $\{Y_i\}\subset\mathfrak{p}$ so that $\omega_G=\sum_i(-\tilde{X}^2_i+\tilde{Y}^2_i)$, where $\tilde{X}_i$ and $\tilde{Y}_i$ are the left invariant vector fields on $G$ associated with $X_i$ and $Y_i$ respectively. We set $\omega'_K:=-2\sum_i\tilde{X}^2_i$. Local coordinates on $G\times K/H$ near $\pi(g,k)$ are given by the inverse of the map
\begin{equation*}
\begin{split}
&(s,t,t_K)=(\ldots,s_i,\ldots,t_i,\ldots,t_{K,i},\ldots)\mapsto\\
&\pi(g\exp(tX+sY),k\exp(t_KX_K)):=\pi(g\exp(\sum_it_iX_i+s_iY_i),k\exp(\sum_it_{K,i}X_{K,i}))
\end{split}
\end{equation*}
for $(s,t,t_K)$ in a neighborhood of the origin in Euclidean space.

Let $C^{\infty}(G\times K/H)$ be the set of smooth functions on $G\times K/H$. We denote by $C^{\infty}(G\times K)^H:=\{h\in C^{\infty}(G\times K)\ ;\ h(gx,kx)=h(g,k)\textrm{ for }x\in K\}$ the set of right-$H$-invariant functions on $G\times K$. For a function $f=f((g,k)H)\in C^{\infty}(G\times K/H)$, we write $\tilde{f}:=f\circ\pi\in C^{\infty}(G\times K)^H$. We define $D_{\omega_G\otimes1}$ and $D_{1\otimes\omega_K}$ in the algebra $\mathbf{D}(G\times K/H)$ of left-$G\times K$-invariant differential operators. We have the orthogonal complement $\mathfrak{h}^{\perp}\subset\mathfrak{g}\oplus\mathfrak{k}$ with respect to the bilinear form $B\oplus B_{\mathfrak{k}}$. Let $\operatorname{Sym}(\mathfrak{h}^{\perp})$ be the symmetric algebra on $\mathfrak{h}^{\perp}$, and $I(\mathfrak{h}^{\perp})\subset \operatorname{Sym}(\mathfrak{h}^{\perp})$ be the $\operatorname{Ad}_{G\times K}(H)$-invariant subalgebra. Note that, making an identification of $\operatorname{Sym}(\mathfrak{g}\oplus\mathfrak{k})$ with the commutative algebra $\mathbb{C}[X_i,Y_i,X_{K,i}]$ of polynomials, we can take a polynomial $Q(X,Y,X_K)$ corresponding to each $Q\in\operatorname{Sym}(\mathfrak{g}\oplus\mathfrak{k})$. For $Q\in\operatorname{Sym}(\mathfrak{h}^{\perp})\subset\operatorname{Sym}(\mathfrak{g}\oplus\mathfrak{k})$, a polynomial $Q(X,Y,X_K)$ representing $Q$ can be taken independently of the directions in $\mathfrak{h}$. Cf. \cite[Chapter II, Lemma 4.7]{Hel00}. It is shown in \cite[Chapter II, Theorem 4.9]{Hel00} that there exists a unique linear bijection of $I(\mathfrak{h}^{\perp})$ onto $\mathbf{D}(G\times K/H)$ such that the operator $D_{\tilde{Q}}\in\mathbf{D}(G\times K/H)$ for $Q\in I(\mathfrak{h}^{\perp})$ is described as
\begin{equation*}
(D_{\tilde{Q}}f)((g,k)H)=(\left.Q(\partial_s,\partial_t,\partial_{t_K})\right|_{(s,t,t_K)=0}\tilde{f})(g\exp(tX+sY),k\exp(t_KX_K))
\end{equation*}
for $f\in C^{\infty}(G\times K/H)$. Thus, $D_{\omega_G\otimes1}$ and $D_{1\otimes\omega_K}$ in $\mathbf{D}(G\times K/H)$ are respectively defined as
\begin{align}
&(D_{\omega_G\otimes1}f)((g,k)H):=(\left.-\partial^2_t+\partial^2_s\right|_{(s,t)=0}\tilde{f})(g\exp(tX+sY),k),\label{SpecialFunctions:eq3.1.1}\\
&(D_{1\otimes\omega_K}f)((g,k)H):=(\left.-2\partial^2_{t_K}\right|_{t_K=0}\tilde{f})(g,k\exp(t_KX_K)),\label{SpecialFunctions:eq3.1.2}
\end{align}
where $\partial^2_{\bullet}:=\sum_i \partial^2_{\bullet_i}$. Note that, under the the linear bijection of $I(\mathfrak{h}^{\perp})$ onto $\mathbf{D}(G\times K/H)$, cf. \cite[Chapter II, Theorem 4.9]{Hel00}, $D_{-\omega_G\otimes1+1\otimes\omega_K}=-D_{\omega_G\otimes1}+D_{1\otimes\omega_K}$.

\subsection{Explicit formulae for $\rho_t^{G_0}$ and $\rho_t^G$}
\subsubsection{A solution to the heat equation}
We define a homogeneous space, denoted by $P$, which is a $G\times K$-space constructed via the Iwasawa projection onto $NA\subset G$ together with a multiplicative cocycle from $G\times NA$ to $K$ associated with $NA$. Here, the Iwasawa projection $\operatorname{Iw}_K$ to $K$ gives rise to a map
\begin{equation*}
\kappa:\ G\times NA\to K,\ (g,na)\mapsto \operatorname{Iw}_K(gna),
\end{equation*}
which satisfies the cocycle identity, for $g_1,g_2\in G$ and $na\in NA$,
\begin{equation*}
\kappa(g_1g_2,na)=\kappa(g_1,g_2na)\kappa(g_2,na)
\end{equation*}
by uniqueness of the Iwasawa decomposition $G=NAK$. Using Lemma~{\upshape\ref{SpecialFunctions:lem3.1.1}}, we choose a representative $[g,k]\in G\times K$ for every element in $G\times K/H$, and set $P:=\{[g,k]\ ;\ g\in G,\ k\in K\}$ with the $G\times K$-action $(g',k')\cdot[g,k]:=[g'g,k'k]$.

\begin{lem}\label{SpecialFunctions:lem3.1.1}
For $g\in G$ and $k\in K$, let $\operatorname{Iw}_{NA}(g):=\operatorname{Iw}_N(g)\operatorname{Iw}_A(g)$, $[g,k]:=\operatorname{Iw}(g,k):=(\operatorname{Iw}_{NA}(gk^{-1}),(\operatorname{Iw}_K(gk^{-1}))^{-1})$. Then for any $(g,k)H\in G\times K/H$, we have $[g,k]H=(g,k)H$.
\end{lem}
\begin{proof}
Take $(gx,kx)\in(g,k)H$ with $x\in K$. Since $\operatorname{Iw}_{NA}(g)=\operatorname{Iw}_{NA}(gk^{-1})$ and $\operatorname{Iw}_K(g)k^{-1}=\operatorname{Iw}_K(gk^{-1})$, the right translation map
\begin{equation*}
(g',k')\mapsto(g'(\operatorname{Iw}_K(g)x)^{-1},k'(\operatorname{Iw}_K(g)x)^{-1})
\end{equation*}
of $G\times K$ gives
\begin{equation*}
\begin{split}
(gx(\operatorname{Iw}_K(g)x)^{-1},kx(\operatorname{Iw}_K(g)x)^{-1})&=(\operatorname{Iw}_{NA}(g),k(\operatorname{Iw}_K(g))^{-1})\\
&=(\operatorname{Iw}_{NA}(gk^{-1}),(\operatorname{Iw}_K(gk^{-1}))^{-1}).
\end{split}
\end{equation*}
Thus, the classes coincide in $G\times K/H$.
\end{proof}

Similarly to the description of the differential operators $\mathbf{D}(G\times K/H)$, we present the differential operators on the quotient space $P$. We denote by $\mathbf{D}(G\times K)$ the algebra of left-$G\times K$-invariant differential operators on the Lie group $G\times K$. Let $C^{\infty}(P)$ be the set of smooth functions on $P$. We write such a function as $\varphi([g,k]):=\varphi(\operatorname{Iw}_{NA}(gk^{-1}),(\operatorname{Iw}_K(gk^{-1}))^{-1})$. Then using the $G\times K$-equivariant diffeomorphism $(g,k)H\mapsto [g,k]$, the algebra of left-$G\times K$-invariant differential operators on $P$ is identified with $\mathbf{D}(G\times K/H)$ on the homogeneous space $G\times K/H$. Hence $D_{\tilde{Q}}\in\mathbf{D}(G\times K/H)$ for $Q\in I(\mathfrak{h}^{\perp})$ acts on $C^{\infty}(P)$ by
\begin{equation*}
\begin{split}
&(D_{\tilde{Q}}\varphi)([g,k])\\
&=(\left.Q(\partial_s,\partial_t,\partial_{t_K})\right|_{(s,t,t_K)=0}\bar{\varphi})(\operatorname{Iw}_{NA}(gk^{-1})\exp(tX+sY),(\operatorname{Iw}_K(gk^{-1}))^{-1}\exp(t_KX_K))
\end{split}
\end{equation*}
for $\varphi\in C^{\infty}(P)$, where we put $\bar{\varphi}:=\varphi\circ\operatorname{Iw}\in C^{\infty}(G\times K)^H$. Thus, corresponding respectively to the two differential operators 
\eqref{SpecialFunctions:eq3.1.1} and \eqref{SpecialFunctions:eq3.1.2} in $\mathbf{D}(G\times K/H)$, 
we define $D_{\omega_G\otimes1}, D_{1\otimes\omega_K}\in\mathbf{D}(P)$ by
\begin{equation*}
\begin{split}
&(D_{\omega_G\otimes1}\varphi)([g,k])=(\left.-\partial^2_t+\partial^2_s\right|_{(s,t)=0}\bar{\varphi})(\operatorname{Iw}_{NA}(gk^{-1})\exp(tX+sY),(\operatorname{Iw}_K(gk^{-1}))^{-1}),\\
&(D_{1\otimes\omega_K}\varphi)([g,k])=(\left.-2\partial^2_{t_K}\right|_{t_K=0}\bar{\varphi})(\operatorname{Iw}_{NA}(gk^{-1}),(\operatorname{Iw}_K(gk^{-1}))^{-1}\exp(t_KX_K))
\end{split}
\end{equation*}
for $\varphi\in C^{\infty}(P)$.

If $g\in G$, let $L_g$ and $R_g$ denote the left and right translations $x\mapsto gx$ and $x\mapsto xg$ respectively. We denote by $\imath:\ G\ni g\mapsto[g,e]\in P$ the diffeomorphism. By definition, it follows that $\imath\circ L_g=L_{(g,e)}\circ\imath$ for $g\in G$. We denote by $\mathbf{D}(G)$ the algebra of left-$G$-invariant differential operators on $G$. $D\in\mathbf{D}(G)$ is called right-$K$-invariant if $D^{R_k}=D$ for all $k\in K$, where we put $D^{R_k}f=(D(f\circ R_k))\circ R_{k^{-1}}$ for $f\in C^{\infty}(G)$. Then Lemma~{\upshape\ref{SpecialFunctions:lem3.1.2}} shows that, for a right-$K$-invariant $D\in\mathbf{D}(G)$, the tangent map $d\imath_e$ at $e\in G$ induces a push-forward operator on $P$.
\begin{lem}\label{SpecialFunctions:lem3.1.2}
Let $Q$ be a polynomial in the symmetric algebra $\operatorname{Sym}(\mathfrak{g})$ on $\mathfrak{g}$, and $D=D_{\tilde{Q}}\in\mathbf{D}(G)$ be a differential operator on $G$ with
\begin{equation*}
Df(g)=(\left.Q(\partial_t,\partial_s)\right|_{(s,t)=0}f)(g\exp(tX+sY))
\end{equation*}
for $f\in C^{\infty}(G)$. Then if $D\in\mathbf{D}(G)$ is right-$K$-invariant, we have
\begin{equation}\label{SpecialFunctions:eq3.1.3}
\begin{split}
&((Df)\circ\imath^{-1})([g,k])\\
&=(\left.Q(\partial_t,\partial_s)\right|_{(s,t)=0}\overline{(f\circ\iota^{-1})})(\operatorname{Iw}_{NA}(gk^{-1})\exp(tX+sY),(\operatorname{Iw}_K(gk^{-1}))^{-1}).
\end{split}
\end{equation}
\end{lem}

\begin{proof}
Since $D^{R_x}=D$ for $x\in K$, $D$ is $\operatorname{Ad}_G(K)$-invariant, hence we have another basis $\{\operatorname{Ad}_G(x)X_i,\operatorname{Ad}_G(x)Y_i\}\subset \mathbf{D}(G)$ for each $x\in K$. Let us analyze the left hand side of the desired equation \eqref{SpecialFunctions:eq3.1.3}:
\begin{equation*}
\begin{split}
&((Df)\circ\imath^{-1})([g,k])=Df(gk^{-1})=((Df)\circ L_{gk^{-1}})(e)=(D(f\circ L_{(gk^{-1})^{-1}}))(e).
\end{split}
\end{equation*}
Applying the $\operatorname{Ad}_G(K)$-invariance of $D$, we get 
\begin{equation*}
\begin{split}
&(D(f\circ L_{(gk^{-1})^{-1}}))(e)=\operatorname{Ad}_G(k(\operatorname{Iw}_K(g))^{-1})D(f\circ L_{(gk^{-1})^{-1}})(e)\\
&=(\left.Q(\partial_t,\partial_s)\right|_{(s,t)=0}(f\circ L_{(gk^{-1})^{-1}}))(\exp(t\operatorname{Ad}_G(k(\operatorname{Iw}_K(g))^{-1})X+s\operatorname{Ad}_G(k(\operatorname{Iw}_K(g))^{-1})Y))\\
&=(\left.Q(\partial_t,\partial_s)\right|_{(s,t)=0}(f\circ L_{(gk^{-1})^{-1}}\circ\imath^{-1}))([\operatorname{Ad}_G(k(\operatorname{Iw}_K(g))^{-1})\exp(tX+sY),e]).
\end{split}
\end{equation*}
Now, using the property that $L_{(gk^{-1})^{-1}}\circ\imath^{-1}=\imath^{-1}\circ L_{((gk^{-1})^{-1},e)}$, we have
\begin{equation*}
\begin{split}
&=(\left.Q(\partial_t,\partial_s)\right|_{(s,t)=0}(f\circ L_{(gk^{-1})^{-1}}\circ\imath^{-1}))([\operatorname{Ad}_G(k(\operatorname{Iw}_K(g))^{-1})\exp(tX+sY),e])\\
&=(\left.Q(\partial_t,\partial_s)\right|_{(s,t)=0}(f\circ\imath^{-1}\circ L_{((gk^{-1})^{-1},e)}))([\operatorname{Ad}_G(k(\operatorname{Iw}_K(g))^{-1})\exp(tX+sY),e])\\
&=(\left.Q(\partial_t,\partial_s)\right|_{(s,t)=0}f\circ\imath^{-1})([\operatorname{Iw}_{NA}(gk^{-1})\exp(tX+sY)\operatorname{Iw}_K(gk^{-1}),e]).
\end{split}
\end{equation*}
Finally, the definition of the map $\iota$ and the identity $\exp(tX+sY)=\operatorname{Iw}_{NA}(\exp(tX+sY))\operatorname{Iw}_K(\exp(tX+sY))$ yields
\begin{equation*}
\begin{split}
&(\left.Q(\partial_t,\partial_s)\right|_{(s,t)=0}f\circ\imath^{-1})([\operatorname{Iw}_{NA}(gk^{-1})\exp(tX+sY)\operatorname{Iw}_K(gk^{-1}),e])\\
&=(\left.Q(\partial_t,\partial_s)\right|_{(s,t)=0}\overline{(f\circ\imath^{-1})})(\operatorname{Iw}_{NA}(gk^{-1})\exp(tX+sY),(\operatorname{Iw}_K(gk^{-1}))^{-1}),
\end{split}
\end{equation*}
which completes the proof.
\end{proof}

We are now ready to deduce the heat kernel for the Laplace--Beltrami operator $\Delta_G$ on $G$. To describe $\Delta_G$ in terms of the differential operators on the manifold $P=\{[g,k]=\operatorname{Iw}(g,k)\ ;\ g\in G,k\in K\}$, we define the push-forward for the diffeomorphism $\imath:\ G\to P$, using Lemma~{\upshape\ref{SpecialFunctions:lem3.1.2}}.
\begin{definition}\label{SpecialFunctions:def3.2.1}
Let $D=D_{\tilde{Q}}\in\mathbf{D}(G)$ be a right-$K$-invariant differential operator for $Q\in\operatorname{Sym}(\mathfrak{g})$. Then the push-forward $\imath_{*}D$ is defined to be
\begin{equation*}
\begin{split}
&\imath_{*}D\varphi([g,e]):=(D(\varphi\circ\imath)\circ\imath^{-1})([g,e])=D(\varphi\circ\imath)(g)\\
&=(\left.Q(\partial_t,\partial_s)\right|_{(s,t)=0}\bar{\varphi})(\operatorname{Iw}_{NA}(g)\exp(tX+sY),(\operatorname{Iw}_K(g))^{-1})
\end{split}
\end{equation*}
for $\varphi\in C^{\infty}(P)$ and $g\in G$.
\end{definition}
\noindent Hence we get differential operators $\imath_{*}\Delta_G$, $\imath_{*}\omega_G$ and $\imath_{*}\omega'_K$ on $P$. We now proceed to prove Proposition~{\upshape\ref{SpecialFunctions:prop3.2.1}} that gives an explicit description for $\imath_{*}\Delta_G$.

\begin{prop}\label{SpecialFunctions:prop3.2.1}
As elements in $\mathbf{D}(P)$, we have $\imath_{*}\Delta_G=D_{-\omega_G\otimes1+1\otimes\omega_K}$.
\end{prop}

\begin{proof}
Let $g=nak=n(g)a(g)k(g)\in NAK$, and $\varphi\in C^{\infty}(P)$. Select a representative in $G\times K/H$ so that its second entry is $e\in K$, and then we shall show that
\begin{equation}\label{SpecialFunctions:eq3.2.1}
-D_{\omega_G\otimes1}\varphi([g,e])+D_{1\otimes\omega_K}\varphi([g,e])=\imath_{*}\Delta_G\varphi([g,e]).
\end{equation}
By definition, the first term of the left hand side of the equation \eqref{SpecialFunctions:eq3.1.1} can be written as
\begin{equation*}
D_{\omega_G\otimes1}\varphi([g,e])=(\left.-\partial^2_t+\partial^2_s\right|_{(s,t)=0}\bar{\varphi})(na\exp(tX+sY),k^{-1})=\imath_{*}\omega_G\varphi(g).
\end{equation*}
For the second term of the left hand side, since $\operatorname{Ad}_K(k)\omega_K=\Delta^{R_k}_K=\omega_K$, we have
\begin{equation*}
\begin{split}
&D_{1\otimes\omega_K}\varphi([g,e])=D_{1\otimes\operatorname{Ad}_K(k)\omega_K}\varphi([g,e])\\
&=(\left.-\partial^2_{t_K}\right|_{t_K=0}\bar{\varphi})(na,k^{-1}\exp(\sum_it_{K,i}\operatorname{Ad}_K(k)X_{K,i}))\\
&=(\left.-\partial^2_{t_K}\right|_{t_K=0}\bar{\varphi})(na,\exp(t_KX_K)k^{-1})\\
&=(\left.-\partial^2_{t_K}\right|_{t_K=0}\bar{\varphi})(\operatorname{Iw}_{NA}(g\exp(-t_KX)),(\operatorname{Iw}_K(g\exp(-t_KX)))^{-1})\\
&=(\left.-\partial^2_{t_K}\right|_{t_K=0}\varphi)([g\exp(-t_KX),e])=(\left.-\partial^2_{t_K}\right|_{t_K=0}\varphi\circ\imath)(g\exp(-t_KX))\\
&=(\left.-\partial^2_{t_K}\right|_{(s,t_K)=0}\varphi\circ\imath)(g\exp(-t_KX-sY))=\omega'_K(\varphi\circ\imath)(g)=\imath_{*}\omega'_K\varphi([g,e]).
\end{split}
\end{equation*}
Hence the equation \eqref{SpecialFunctions:eq3.2.1} follows from
\begin{equation*}
-D_{\omega_G\otimes1}\varphi([g,e])+D_{1\otimes\omega_K}\varphi([g,e])=-\imath_{*}\omega_G\varphi([g,e])+\imath_{*}\omega'_K\varphi([g,e])=\imath_{*}\Delta_G\varphi([g,e]).
\end{equation*}
\end{proof}
Given functions $f:\ G\to\mathbb{C}$ and $h:\ K\to\mathbb{C}$, we define their associated function on $P$ with the product $fh$ on $G\times K$. That is, for $[g,k]\in P$, we define
\begin{equation*}
[f\times h]([g,k]):=f(\operatorname{Iw}_{NA}(g))h(k(\operatorname{Iw}_K(g))^{-1}).
\end{equation*}

\begin{proof}[Proof of Theorem~{\upshape\ref{SpecialFunctions:thm3.0.1}}]
It will suffice to find a solution to the problem \eqref{SpecialFunctions:eq2.1.1} by the uniqueness of the heat kernel $\rho^G_t$. We construct such a solution using the differential operator $-\omega_G\otimes1+1\otimes\omega_K$ on $G\times K$. Under the diffeomorphism $\imath:\ G\to P$, it follows from Proposition~{\upshape\ref{SpecialFunctions:prop3.2.1}} that $\rho^G_t$ on $G$ is induced by a certain solution for the differential operator $D_{-\omega_G\otimes1+1\otimes\omega_K}$ on $P$. Moreover, by using the definition of $D_{-\omega_G\otimes1+1\otimes\omega_K}$, a solution for $D_{-\omega_G\otimes1+1\otimes\omega_K}$ in $C^{\infty}(P)$ corresponds to a solution in $C^{\infty}(G\times K)^H$ for $-\omega_G\otimes1+1\otimes\omega_K$.

With respect to the Cartan decomposition, every element $h\in G$ can be uniquely written as $h=\sqrt{h(\sigma(h))^{-1}}k_h=(\operatorname{Ad}_G(\operatorname{Crt}_K(h))\operatorname{Crt}_{\overline{A^{+}}}(h))k_h\in(\operatorname{Ad}_G(K)\overline{A^{+}})K$, where $k_h=\operatorname{Crt}_K(h)\operatorname{Crt}'_K(h)$. Recall that each spherical function $\varphi_{\lambda}$ in the equation \eqref{SpecialFunctions:eq2.1.2} is right-$K$-invariant. Define
\begin{equation*}
(g^G_{t/2})_{\operatorname{Ad}_G(K)\overline{A^{+}}}(h):=g^G_{t/2}(\sqrt{h(\sigma(h))^{-1}}),
\end{equation*}
and set
\begin{equation}\label{SpecialFunctions:eq3.2.2}
\tilde{p}_t(h,k):=(g^G_{t/2})_{\operatorname{Ad}_G(K)\overline{A^{+}}}(h)\rho^K_t(kk^{-1}_h).
\end{equation}
Here, for a right-$K$-invariant function $F$ on $G$, the definition
\begin{equation*}
F_{\operatorname{Ad}_G(K)\overline{A^{+}}}(h):=F(\sqrt{h(\sigma(h))^{-1}})
\end{equation*}
for $h=\sqrt{h(\sigma(h))^{-1}}k_h\in G$ gives a well-defined function on $G$, and $F_{\operatorname{Ad}_G(K)\overline{A^{+}}}$ may be regarded as the composition of $F$ and the projection $G\to\operatorname{Ad}_G(K)\overline{A^{+}}$. By definition, we get
\begin{equation*}
\begin{split}
\tilde{p}_t(hx,kx)&=(g^G_{t/2})_{\operatorname{Ad}_G(K)\overline{A^{+}}}(hx)\rho^K_t(kx(k_{hx})^{-1})\\
&=(g^G_{t/2})_{\operatorname{Ad}_G(K)\overline{A^{+}}}(h)\rho^K_t(kx(k_hx)^{-1})=\tilde{p}_t(h,k)
\end{split}
\end{equation*}
for $x\in K$. Hence, combining this with $g^G_{t/2}\in C^{\infty}(G)^K$ and $\rho^K_t\in C^{\infty}(K)$, we obtain $\tilde{p}_t\in C^{\infty}(G\times K)^H$. Using Lemma~{\upshape\ref{Chap2:lem2.1.1}}, which is proved at the end of this subsubsection, we obtain the equation
\begin{equation*}
(\partial_t-\omega_G\otimes1+1\otimes\omega_K)\tilde{p}_t(h,k)=0.
\end{equation*}

\begin{lem}\label{Chap2:lem2.1.1}
$\tilde{p}_t$ satisfies the heat equation associated with $-\omega_G\otimes1+1\otimes\omega_K$ on $G\times K/H$.
\end{lem}

Recall that the integrand of a spherical function $\varphi_{\lambda}(h)$ for $\lambda\in i\mathfrak{a}^{*}$ depends only on the $A$-component of $h\in G$ in the Iwasawa decomposition $G=NAK$. By the equation \eqref{SpecialFunctions:eq2.1.2}, we have
\begin{equation*}
\begin{split}
(g^G_{t/2})_{\operatorname{Ad}_G(K)\overline{A^{+}}}(h)&=\frac{1}{|W|}\int_K\int_{i\mathfrak{a}^{*}}e^{(\langle\lambda,\lambda\rangle-\langle\rho,\rho\rangle)\frac{t}{2}}e^{(\lambda+\rho)\log\operatorname{Iw}_A(\sqrt{xh(\sigma(xh))^{-1}})}\frac{d\lambda}{|\mathbf{c}(\lambda)|^2}dx\\
&=\frac{1}{|W|}\int_K\int_{i\mathfrak{a}^{*}}e^{(\langle\lambda,\lambda\rangle-\langle\rho,\rho\rangle)\frac{t}{2}}e^{(\lambda+\rho)\log\operatorname{Iw}_A(x\sqrt{h(\sigma(h))^{-1}}x^{-1})}\frac{d\lambda}{|\mathbf{c}(\lambda)|^2}dx.
\end{split}
\end{equation*}
Here, the first equality follows from Fubini's theorem. Under the $G\times K$-equivariant diffeomorphism
\begin{equation*}
G\times K/H\to P,\ (g,k)H\mapsto [g,k]=\operatorname{Iw}(g,k)=(\operatorname{Iw}_{NA}(gk^{-1}),(\operatorname{Iw}_K(gk^{-1}))^{-1}),
\end{equation*}
we pick the representative
\begin{equation*}
\begin{split}
&(\operatorname{Iw}_{NA}(\sqrt{xh(\sigma(xh))^{-1}}k_hk^{-1}),(\operatorname{Iw}_K(\sqrt{xh(\sigma(xh))^{-1}}k_hk^{-1}))^{-1})\\
&=(\operatorname{Iw}_{NA}(x\sqrt{h(\sigma(h))^{-1}}x^{-1}),(\operatorname{Iw}_K(x\sqrt{h(\sigma(h))^{-1}}x^{-1}k_hk^{-1}))^{-1})
\end{split}
\end{equation*}
of $(\sqrt{xh(\sigma(xh))^{-1}},kk^{-1}_h)H\in G\times K/H$ for each $x\in K$, with respect to the function \eqref{SpecialFunctions:eq3.2.2} in $C^{\infty}(G\times K)^H$. From this, we get a solution
\begin{equation}\label{SpecialFunctions:eq3.2.3}
\begin{split}
[h,k]&\mapsto [g^G_{t/2}\times\rho^K_t]([\operatorname{Iw}_{NA}(\sqrt{h(\sigma(h))^{-1}}),(\operatorname{Iw}_K(\sqrt{h(\sigma(h))^{-1}}k_hk^{-1}))^{-1}])\\
&=\frac{1}{|W|}\int_K\left(\int_{i\mathfrak{a}^{*}}e^{(\langle\lambda,\lambda\rangle-\langle\rho,\rho\rangle)\frac{t}{2}}e^{(\lambda+\rho)\log\operatorname{Iw}_A(x\sqrt{h(\sigma(h))^{-1}}x^{-1})}\frac{d\lambda}{|\mathbf{c}(\lambda)|^2}\right)\\
&\hspace{20mm}\times\rho^K_t((\operatorname{Iw}_K(x\sqrt{h(\sigma(h))^{-1}}x^{-1}k_hk^{-1}))^{-1})dx
\end{split}
\end{equation}
of the heat equation associated with $D_{-\omega_G\otimes1+1\otimes\omega_K}$ on $P$. Then Proposition~{\upshape\ref{SpecialFunctions:prop3.2.1}} implies that, under the equivariant diffeomorphism from $G$-space $G$ to $G\times\{e\}$-space $P$, we can identify the solution \eqref{SpecialFunctions:eq3.2.2} with the solution $p^G_t$, defined by
\begin{equation}\label{SpecialFunctions:eq3.2.4}
\begin{split}
p^G_t(h)&:=[g^G_{t/2}\times\rho^K_t]([\operatorname{Iw}_{NA}(\sqrt{h(\sigma(h))^{-1}}),(\operatorname{Iw}_K(\sqrt{h(\sigma(h))^{-1}}k_h))^{-1}])\\
&=\frac{1}{|W|}\int_K\left(\int_{i\mathfrak{a}^{*}}e^{(\langle\lambda,\lambda\rangle-\langle\rho,\rho\rangle)\frac{t}{2}}e^{(\lambda+\rho)\log\operatorname{Iw}_A(x\sqrt{h(\sigma(h))^{-1}}x^{-1})}\frac{d\lambda}{|\mathbf{c}(\lambda)|^2}\right)\\
&\hspace{20mm}\times\rho^K_t((\operatorname{Iw}_K(x\sqrt{h(\sigma(h))^{-1}}x^{-1}k_h))^{-1})dx,
\end{split}
\end{equation}
of the heat equation of the problem \eqref{SpecialFunctions:eq2.1.1}. Furthermore, by Lemma~{\upshape\ref{SpecialFunctions:lem3.2.1}}, the pullback $p^G_t$ by $\imath$ of the function \eqref{SpecialFunctions:eq3.2.2} satisfies the initial condition of the problem \eqref{SpecialFunctions:eq2.1.1}. Therefore, the pullback $p^G_t$ gives the desired expression \eqref{SpecialFunctions:eq3.0.1}.
\end{proof}
By a similar argument, we obtain the formula \eqref{SpecialFunctions:eq3.0.1} for $\rho_t^{G_0}$.\\

To show Lemma~{\upshape\ref{Chap2:lem2.1.1}}, we need some notation for a decomposition of $L^2(G\times K)$. Using the Peter--Weyl theorem for $K$, we get
\begin{equation*}
L^2(K,dk)\simeq \bigoplus_{\tau\in\hat{K}} V_{\tau}\otimes V^{*}_{\tau},
\end{equation*}
where $(\tau^{\wedge},V^{*}_{\tau})$ denotes the contragredient representation of $(\tau,V_{\tau})$ with respect to a Hermitian form $\langle\cdot,\cdot\rangle_{\tau}$ that is anti-linear on the second coordinate. For each $\tau\in\hat{K}$, take an orthonormal basis $\{v_{\tau,1},\ldots,v_{\tau,\dim V_{\tau}}\}$ of $V_{\tau}$. We denote by $\{v^{*}_{\tau,1},\ldots,v^{*}_{\tau,\dim V_{\tau}}\}$ the dual basis of $V^{*}_{\tau}$. Then $L^2(K)$ admits an orthonormal basis consisting of the functions
\begin{equation*}
u^{\tau}_{i,j}(k):=\langle \tau(k)v_{\tau,i},v^{*}_{\tau,j}\rangle_{\tau},
\end{equation*}
where $k\in K$, $\tau\in\hat{K}$ and $i,j=1,\ldots,\dim V_{\tau}$. For each $\tau\in\hat{K}$, let $\operatorname{proj}_{\tau}:\ L^2(K,dk)\to V_{\tau}\otimes V^{*}_{\tau}$ be an orthogonal projection, and let
\begin{equation*}
Q_{\tau}:=\left(\int_K\mathcal{R}(k)\otimes\operatorname{id}_{V_{\tau}}\otimes\tau^{\wedge}(k)dk\right)\circ(\operatorname{id}_{L^2(G)}\otimes\operatorname{proj}_{\tau})
\end{equation*}
be the orthogonal projection of $L^2(G\times K)\simeq L^2(G,dg)\otimes L^2(K,dk)$ onto the $H$-invariant subspace $(L^2(G,dg)\otimes V_{\tau}\otimes V^{*}_{\tau})^H$, where $\mathcal{R}$ denotes the right regular representation of $G$ on $L^2(G,dg)$. For $\tau\in\hat{K}$, we have
\begin{equation*}
Q_{\tau}(-\omega_G\otimes1+1\otimes\omega_K)Q_{\tau}=Q_{\tau}(-\omega_G\otimes1)Q_{\tau}+\lambda_{\tau}(\operatorname{id}_{L^2(G)}\otimes\operatorname{proj}_{\tau})
\end{equation*}
on $L^2(G\times K)$.

We recall the measure on $\operatorname{Ad}_G(K)\overline{A^{+}}\mathbf{o}$ induced by the Haar measure $dg$ on $G$. Here we denote by $\mathbf{o}$ the origin $eK$ in the symmetric space $G/K$. Let
\begin{equation*}
\delta(\exp H):=
\prod_{\alpha\in\Sigma^{+}}\left(\frac{2\sinh\alpha(H)}{\langle\alpha,\rho\rangle}\right)^{m_{\alpha}}
\end{equation*}
be the Jacobian determinant on the Weyl chamber $H\in\mathfrak{a}^{+}$. We denote by $Z_K(\mathfrak{a})$ the centralizer of $\mathfrak{a}$ in $K$. It is shown in \cite[Chapter I, Theorem 5.8]{Hel00} that the symmetric space $\operatorname{Ad}_G(K)\overline{A^{+}}\mathbf{o}\approx G/K$ admits the measure $d\mu$ given by the integration formula
\begin{equation}\label{Chap2:eq2.1.1}
\int_{\operatorname{Ad}_G(K)\overline{A^{+}}\mathbf{o}}f(kak^{-1}\mathbf{o})d\mu(kak^{-1}\mathbf{o})=C\int_{K/Z_K(\mathfrak{a})}\left(\int_{A^{+}}f(kak^{-1}\mathbf{o})\delta(a)da\right)d(kZ_K(\mathfrak{a}))
\end{equation}
for $f\in C^{\infty}_c(\operatorname{Ad}_G(K)\overline{A^{+}}\mathbf{o})$, where $C$ is a normalization constant dependent on the choice of measures, and $d(kZ_K(\mathfrak{a}))$ denotes a $K$-invariant measure on the quotient space $K/Z_K(\mathfrak{a})$ that is normalized by $\int d(kZ_K(\mathfrak{a}))=1$.

Using the Cartan decomposition $G=K\overline{A^{+}}K$, we get a diffeomorphism $\iota:\ \operatorname{Ad}_G(K)\overline{A^{+}}\times K\to G,\ (kak^{-1},k')\mapsto kak^{-1}k'$.
\begin{lem}\label{Chap2:lem2.1.2}
The map
\begin{equation}\label{Chap2:eq2.1.2}
\begin{split}
&L^2(\operatorname{Ad}_G(K)\overline{A^{+}}\mathbf{o},d\mu)\otimes V_{\tau}\otimes V^{*}_{\tau}\to (L^2(G,dg)\otimes V_{\tau}\otimes V^{*}_{\tau})^H,\\
&\varphi\otimes u^{\tau}_{ij}\mapsto ((\iota(kak^{-1},k'),x)\mapsto \varphi(kak^{-1})\langle \tau(x)v_{\tau,i},\tau^{\wedge}(k')v^{*}_{\tau,j}\rangle_{\tau})
\end{split}
\end{equation}
for $i,j=1,\ldots,\dim V_{\tau}$ induces an isomorphism of left-$G\times K$ regular representations
\begin{equation}\label{Chap2:eq2.1.3}
L^2(\operatorname{Ad}_G(K)\overline{A^{+}}\mathbf{o},d\mu)\otimes V_{\tau}\otimes V^{*}_{\tau}\simeq (L^2(G,dg)\otimes V_{\tau}\otimes V^{*}_{\tau})^H.
\end{equation}
\end{lem}
\begin{proof}
For each function $f \in(L^2(G,dg)\otimes V_{\tau}\otimes V^{*}_{\tau})^H$, we obtain the expansion
\begin{equation*}
f(g,x)=\sum^{\dim V_{\tau}}_{i,j=1}f_{ij}(g)u^{\tau}_{i,j}(x)
\end{equation*}
for $g\in G$ and $x\in K$, where $f_{ij}\in L^2(G,dg)$. Then the inverse of the map \eqref{Chap2:eq2.1.2} is given by
\begin{equation*}
\begin{split}
(L^2(G,dg)\otimes V_{\tau}\otimes V^{*}_{\tau})^H\ni f\mapsto \sum^{\dim V_{\tau}}_{i,j=1}f_{ij}(\iota(\bullet,e))\otimes u^{\tau}_{i,j}\in L^2(\operatorname{Ad}_G(K)\overline{A^{+}}\mathbf{o},d\mu)\otimes V_{\tau}\otimes V^{*}_{\tau}.
\end{split}
\end{equation*}
\end{proof}
\noindent By Lemma~{\upshape\ref{Chap2:lem2.1.2}}, we obtain
\begin{equation}\label{Chap2:eq2.1.4}
\begin{split}
&L^2(G\times K/H)\simeq L^2(G\times K)^H\simeq \left(L^2(G,dg)\otimes \bigoplus_{\tau\in\hat{K}} V_{\tau}\otimes V^{*}_{\tau}\right)^H\\
&\simeq\bigoplus_{\tau\in\hat{K}}(L^2(G,dg)\otimes V_{\tau}\otimes V^{*}_{\tau})^H\simeq \bigoplus_{\tau\in\hat{K}}L^2(\operatorname{Ad}_G(K)\overline{A^{+}}\mathbf{o},d\mu)\otimes V_{\tau}\otimes V^{*}_{\tau}
\end{split}
\end{equation}
as left-$G\times K$ regular representations, where the third isomorphism follows from the fact that Hilbert space completion commutes with the tensor product and the orthogonal direct sum $\bigoplus_{\tau\in\hat{K}}V_{\tau}\otimes V^{*}_{\tau}$. Here we regard the direct sum $\bigoplus_{\tau\in\hat{K}}V_{\tau}\otimes V^{*}_{\tau}$ as the completion of all finite direct sums.

\begin{proof}[Proof of Lemma~{\upshape\ref{Chap2:lem2.1.1}}]
Take $\tau\in\hat{K}$. Denote by $\tilde{p}^{\tau}_t$ the image of
\begin{equation*}
(kak^{-1}\mathbf{o},x)\mapsto e^{-\lambda_{\tau}t}g^G_{t/2}(kak^{-1})\chi_{\tau}(x)=e^{-\lambda_{\tau}t}g^G_{t/2}(kak^{-1})\sum^{\dim V_{\tau}}_{i=1}u^{\tau}_{i,i}(x)
\end{equation*}
in $L^2(\operatorname{Ad}_G(K)\overline{A^{+}}\mathbf{o},d\mu)\otimes V_{\tau}\otimes V^{*}_{\tau}$ under the isomorphism \eqref{Chap2:eq2.1.3}. Since we have, for $i=1,\ldots,\dim V_{\tau}$,
\begin{equation*}
\langle \tau(x)v_{\tau,i},\tau^{\wedge}(k)v^{*}_{\tau,i}\rangle_{\tau}=\langle \overline{\tau(k^{-1})}\tau(x)v_{\tau,i},v^{*}_{\tau,i}\rangle_{\tau}=\langle \tau(kx)v_{\tau,i},v^{*}_{\tau,i}\rangle_{\tau},
\end{equation*}
the function $\tilde{p}^{\tau}_t\in (L^2(G,dg)\otimes V_{\tau}\otimes V^{*}_{\tau})^H$ is expressed as
\begin{equation*}
\tilde{p}^{\tau}_t(h,k)=e^{-\lambda_{\tau}t}(g^G_{t/2})_{\operatorname{Ad}_G(K)\overline{A^{+}}}(h)\chi_{\tau}(kk^{-1}_h).
\end{equation*}
Hence it follows that
\begin{equation*}
\begin{split}
&Q_{\tau}(\partial_t-\omega_G\otimes1+1\otimes\omega_K)Q_{\tau}\tilde{p}^{\tau}_t(h,k)\\
&=(\partial_t-Q_{\tau}(\omega_G\otimes1)Q_{\tau}+\lambda_{\tau}(\operatorname{id}_{L^2(G)}\otimes\operatorname{proj}_{\tau}))\tilde{p}^{\tau}_t(h,k)\\
&=(\partial_t-Q_{\tau}(\omega_G\otimes1)Q_{\tau}+\lambda_{\tau})\tilde{p}^{\tau}_t(h,k)\\
&=\chi_{\tau}(kk^{-1}_h)(\partial_t(e^{-\lambda_{\tau}t}(g^G_{t/2})_{\operatorname{Ad}_G(K)\overline{A^{+}}}(h))\\
&-e^{-\lambda_{\tau}t}(\partial_tg^G_{t/2})_{\operatorname{Ad}_G(K)\overline{A^{+}}}(h)+\lambda_{\tau}e^{-\lambda_{\tau}t}(g^G_{t/2})_{\operatorname{Ad}_G(K)\overline{A^{+}}}(h)))\\
&=0.
\end{split}
\end{equation*}
Since we observe
\begin{equation*}
\tilde{p}_t=\sum_{\tau\in \hat{K}}Q_{\tau}\tilde{p}_t=\sum_{\tau\in \hat{K}}(\dim V_{\tau})\tilde{p}^{\tau}_t,
\end{equation*}
we obtain
\begin{equation*}
\begin{split}
&(\partial_t-\omega_G\otimes1+1\otimes\omega_K)\tilde{p}_t(h,k)\\
&=\left(\sum_{\tau\in\hat{K}}Q_{\tau}(\partial_t-\omega_G\otimes1+1\otimes\omega_K)Q_{\tau}\right)\sum_{\tau\in \hat{K}}(\dim V_{\tau})\tilde{p}^{\tau}_t(h,k)\\
&=\sum_{\tau\in\hat{K}}(\dim V_{\tau}) Q_{\tau}(\partial_t-\omega_G\otimes1+1\otimes\omega_K)Q_{\tau}\tilde{p}^{\tau}_t(h,k)=0,
\end{split}
\end{equation*}
which finishes the proof of Lemma~{\upshape\ref{Chap2:lem2.1.1}}.
\end{proof}

\subsubsection{The initial condition}
This subsubsection completes the proof of Theorem~{\upshape\ref{SpecialFunctions:thm3.0.1}} by verifying that the solution $p^G_t$ of the heat equation \eqref{SpecialFunctions:eq2.1.1}, given by the definition \eqref{SpecialFunctions:eq3.2.4}, satisfies the initial condition of the Cauchy problem on $G$ as stated in Lemma~{\upshape\ref{SpecialFunctions:lem3.2.1}}.

\begin{lem}\label{SpecialFunctions:lem3.2.1}
For $f\in L^2(G)$, we have $\|p^G_t*f-f\|_{L^2}\to0$ as $t\to0+$.
\end{lem}

\begin{proof}
Take $\varepsilon>0$. It is sufficient to assume $f\in C^{\infty}_c(G)$. Let $\Omega:=\operatorname{supp}(f)$ and $\Omega_K:=K\Omega K$ be compact sets. By the anti-Iwasawa decomposition $G=KAN$, the equation \eqref{SpecialFunctions:eq3.2.4} implies that
\begin{equation*}
\begin{split}
&p^G_t*f(h)=\int_K\int_A\int_N p^G_t((kan)^{-1})f(kanh)dkdadn\\
&=\int_K\int_A\int_N[g^G_{t/2}\times\rho^K_t]([(an)^{-1},k])f(kanh)dkdadn\\
&=\int_A\int_N[g^G_{t/2}\times\int_K\rho^K_t(\cdot)f(\cdot anh)dk]([(an)^{-1},k])dadn
\end{split}
\end{equation*}
for $h\in G$. Since $\rho^K_t(k)=\rho^K_t(k^{-1})$, we may also write
\begin{equation*}
p^G_t*f(h)=\int_A\int_N g^G_{t/2}((an)^{-1})(\rho^K_t*f)(anh)dadn=(g^{NA}_{t/2}*(\rho^K_t*f))(h),
\end{equation*}
where $g^{NA}_{t/2}$ denotes the restriction of $g^G_{t/2}$ to $NA\subset G$, and we define
\begin{equation*}
(\rho^K_t*f)(h):=\int_K\rho^K_t(k^{-1})f(kh)dk.
\end{equation*}
Note that, due to the definition of $\operatorname{Iw}(g,k)$ in $P$ for $g\in G$ and $k\in K$, the expression for the convolution $\rho^K_t*f$ with respect to $K$ is inconsistent with that of the convolution product on $G$. We have $\rho^K_t*f=0$ on the complement $\Omega^c_K\subset G$. Putting $f^K_t:=\rho^K_t*f$, we get
\begin{equation*}
\begin{split}
&\|p^G_t*f-f\|_{L^2}=\|g^{NA}_{t/2}*f^K_t-f\|_{L^2}=\left(\int_{\Omega_K}+\int_{\Omega^c_K}|g^{NA}_{t/2}*f^K_t(h)-f(h)|^2dh\right)^{\frac{1}{2}}.
\end{split}
\end{equation*}
Now, we deduce Lemma~{\upshape\ref{SpecialFunctions:lem3.2.1}} from the following convergence results stated in Lemma~{\upshape\ref{SpecialFunctions:lem3.2.2}}, which will be proved later. By Lemma~{\upshape\ref{SpecialFunctions:lem3.2.2}} (ii), there is a $t_{\Omega_K}>0$ such that
\begin{equation*}
\sup_{h\in\Omega_K}|g^{NA}_{t/2}*f^K_t(h)-f(h)|<\varepsilon
\end{equation*}
for $0<t<t_{\Omega_K}$. From Lemma~{\upshape\ref{SpecialFunctions:lem3.2.2}} (iii) it follows that there exists a $t_{\Omega^c_K}>0$ such that, for $0<t<t_{\Omega^c_K}$,
\begin{equation*}
\int_{\Omega^c_K}|g^{NA}_{t/2}*f^K_t(h)-f(h)|^2dh<\varepsilon^2.
\end{equation*}
Thus we conclude that, for $0<t<\min\{t_{\Omega_K},t_{\Omega^c_K}\}$,
\begin{equation*}
\begin{split}
&\|p^G_t*f-f\|_{L^2}\\
&=\left(\int_{\Omega_K}|g^{NA}_{t/2}*f^K_t(h)-f(h)|^2dh+\int_{\Omega^c_K}|g^{NA}_{t/2}*f^K_t(h)-f(h)|^2dh\right)^{\frac{1}{2}}\\
&\leq\left((\sup_{h\in\Omega_K}|g^{NA}_{t/2}*f^K_t(h)-f(h)|)^2\int_{\Omega_K}dh+\int_{\Omega^c_K}|g^{NA}_{t/2}*f^K_t(h)-f(h)|^2dh\right)^{\frac{1}{2}}\\
&<\varepsilon\sqrt{\int_{\Omega_K}dh+1}.
\end{split}
\end{equation*}
\end{proof}

\begin{lem}\label{SpecialFunctions:lem3.2.2}
Let $M:=\max_{x\in G}\{|f(x)|\}<\infty$.
\begin{enumerate}
\item We have the uniform convergence
\begin{equation}\label{SpecialFunctions:eq3.2.5}
f^K_t\to f\textrm{ on }G\textrm{ as }t\to0+.\end{equation}
\item $g^{NA}_{t/2}*f^K_t$ converges uniformly to $f$ on $\Omega_K$ as $t\to0+$.
\item The integral
\begin{equation}\label{SpecialFunctions:eq3.2.6}
\int_{\Omega^c_K}|g^{NA}_{t/2}*f^K_t(h)-f(h)|^2dh
\end{equation}
converges to $0$ as $t\to0+$.
\end{enumerate}
\end{lem}

\begin{proof}[Proof of Lemma~{\upshape\ref{SpecialFunctions:lem3.2.2}}]
We have $\operatorname{supp}(f^K_t)\subset \Omega_K$ for all $t>0$.\\
\noindent (i) Since $f^K_t(h)=0=f(h)$ for $h\in \Omega^c_K$ and $t>0$, we have $\sup_{h\in G}|f^K_t(h)-f(h)|=\sup_{h\in\Omega_K}|f^K_t(h)-f(h)|$ for $t>0$. Hence it suffices to show that $f^K_t$ converges uniformly to $f$ as $t\to0+$ on $\Omega_K$. We denote by $C(K)$ the set of continuous functions on $K$ endowed with the supremum norm $\|\cdot\|_{\textrm{sup}}$. Let $e^{t\Delta_K}$ be the heat semigroup on $L^2(K)$, where $\Delta_K=-\omega_K$. We have $e^{t\Delta_K}F(k)=\int_K\rho^K_t(k'^{-1})F(k'k)dk'$ for $F\in C(K)$. We define the continuous function $f_x(k):=f(kx)$ on $K$ for each $x\in\Omega_K$, and let $\mathscr{F}:=\{f_x\ ;\ x\in\Omega_K\}\subset C(K)$. Then one finds that $\mathscr{F}$ is equicontinuous and pointwise precompact. By the Ascoli--Arzel\`a theorem, the closure $\bar{\mathscr{F}}$ with respect to $\|\cdot\|_{\textrm{sup}}$ is compact in $C(K)$. Select a finite subset $\{F_1,\ldots,F_L\}\subset\bar{\mathscr{F}}$ so that $\bigcup^L_{l=1}B_{\varepsilon}(F_l)$ is an open covering of $\bar{\mathscr{F}}$, where $B_{\varepsilon}(F_l)$ denotes the open ball of radius $\varepsilon$ at $F_l$ for every $1\leq l\leq L$. Pick $\tau=\tau_{\varepsilon,F_1,\ldots,F_L}>0$ so that $\|e^{t\Delta_K}F_l-F_l\|_{\textrm{sup}}<\varepsilon$ for all $l$ if $t<\tau$. Then given $F\in\bar{\mathscr{F}}$, take $F_l$ with $\|F-F_l\|_{\textrm{sup}}<\varepsilon$, and estimate, using the inequality $\inf f\leq e^{t\Delta_K}f(h)\leq \sup f$ for $f\in C(K)$ and $h\in\Omega_K$ (cf. \cite[Theorem 7.16]{Gri09}), 
\begin{equation*}
\begin{split}
\|e^{t\Delta_K}F-F\|_{\textrm{sup}}&\leq \|e^{t\Delta_K}(F-F_l)\|_{\textrm{sup}}+\|F-F_l\|_{\textrm{sup}}+\|e^{t\Delta_K}F_l-F_l\|_{\textrm{sup}}\\
&=2\|F-F_l\|_{\textrm{sup}}+\|e^{t\Delta_K}F_l-F_l\|_{\textrm{sup}}<2\varepsilon+\varepsilon=3\varepsilon.
\end{split}
\end{equation*}
for $t<\tau$. Hence $\sup_{F\in\bar{\mathscr{F}}}\|e^{t\Delta_K}F-F\|_{\textrm{sup}}\to0$ as $t\to0+$.

Given $h\in\Omega_K$, since $f^K_t(h)=\int_K\rho^K_t(k^{-1})f_h(k)dk=e^{t\Delta_K}f_h(e)$, we get the inequality
\begin{equation*}
|f^K_t(h)-f(h)|=\|e^{t\Delta_K}f_h(e)-f_h(e)\|\leq\|e^{t\Delta_K}f_h-f_h\|_{\textrm{sup}}.
\end{equation*}
Thus we have
\begin{equation*}
\sup_{h\in\Omega_K}|f^K_t(h)-f(h)|\leq\sup_{h\in\Omega_K}\|e^{t\Delta_K}f_h-f_h\|_{\textrm{sup}}\leq\sup_{F\in\bar{\mathscr{F}}}\|e^{t\Delta_K}F-F\|_{\textrm{sup}}\to0
\end{equation*}
as $t\to0+$. Therefore we obtain the uniform convergence \eqref{SpecialFunctions:eq3.2.5} on $G$.

\noindent (ii) Since $f$ is continuous on the compact set $\Omega_K$, $f$ is uniformly continuous, i.e. there exists a neighborhood $e\in U\subset G$ such that $|f(gh)-f(h)|<\varepsilon$ for all $h\in\Omega_K$ whenever $g\in U$. Take $\delta_U>0$ so that $\{b\in NA\ ;\ \|b\|<\delta_U\}\ni e$ is contained in $U\cap NA$ and relatively open in $NA$. Note that $\|b\|=d_G(e,b)$ for all $b\in NA=AN$. By Lemma~{\upshape\ref{SpecialFunctions:lem2.1.1}} (i) and (ii), we have for $t>0$
\begin{equation*}
1=\int_K\int_A\int_Ng^G_{t/2}(kan)dkdadn=\int_A\int_Ng^G_{t/2}(an)dadn=\int_A\int_Ng^G_{t/2}((an)^{-1})dadn,
\end{equation*}
and there exists $t_{\delta_U}>0$ such that
\begin{equation*}
\varepsilon>\int_{\{x\in G\ ;\ \|x\|\geq\delta_U\}}g^G_{t/2}(x)dx=\int_{\|an\|\geq\delta_U}g^G_{t/2}((an)^{-1})dadn
\end{equation*}
if $t<t_{\delta_U}$. For $h\in\Omega_K$, we reduce an estimate of the integral
\begin{equation*}
g^{NA}_{t/2}*f^K_t(h)-f(h)
\end{equation*}
to an estimate of the three integrals \eqref{Chap2:eq2.2.1}--\eqref{Chap2:eq2.2.3}
\begin{align}
&\int_A\int_Ng^G_{t/2}((an)^{-1})|f^K_t(anh)-f(anh)|dadn,\label{Chap2:eq2.2.1}\\
&\int_{\|an\|<\delta_U}g^G_{t/2}((an)^{-1})|f(anh)-f(h)|dadn,\label{Chap2:eq2.2.2}\\
&\int_{\|an\|\geq\delta_U}g^G_{t/2}((an)^{-1})|f(anh)-f(h)|dadn\label{Chap2:eq2.2.3},
\end{align}
by means of the inequality
\begin{equation}\label{SpecialFunctions:eq3.2.7}
\begin{split}
&|g^{NA}_{t/2}*f^K_t(h)-f(h)|=\left|\int_A\int_N g^G_{t/2}((an)^{-1})(f^K_t(anh)-f(h))dadn\right|\\
&\leq \int_A\int_N g^G_{t/2}((an)^{-1})|f^K_t(anh)-f(h)|dadn\\
&\leq\int_A\int_Ng^G_{t/2}((an)^{-1})|f^K_t(anh)-f(anh)|dadn\\
&+\int_A\int_Ng^G_{t/2}((an)^{-1})|f(anh)-f(h)|dadn\\
&=\int_A\int_Ng^G_{t/2}((an)^{-1})|f^K_t(anh)-f(anh)|dadn\\
&+\int_{\|an\|<\delta_U}g^G_{t/2}((an)^{-1})|f(anh)-f(h)|dadn\\
&+\int_{\|an\|\geq\delta_U}g^G_{t/2}((an)^{-1})|f(anh)-f(h)|dadn.
\end{split}
\end{equation}

Using the convergence \eqref{SpecialFunctions:eq3.2.5}, the integral \eqref{Chap2:eq2.2.1} is bounded by
\begin{equation*}
\begin{split}
&\sup_{h\in\Omega_K}\int_A\int_Ng^G_{t/2}((an)^{-1})|f^K_t(anh)-f(anh)|dadn\\
&\leq\sup_{h\in\Omega_K}\sup_{b\in AN}|f^K_t(bh)-f(bh)|\int_A\int_Ng^G_{t/2}((an)^{-1})dadn\\
&=\sup_{h\in\Omega_K}\sup_{b\in AN}|f^K_t(bh)-f(bh)|\leq\sup_{h\in G}|f^K_t(h)-f(h)|\to0
\end{split}
\end{equation*}
as $t\to0+$. For $t<t_{\delta_U}$ we have an estimate of the absolute values of the integrals \eqref{Chap2:eq2.2.2} and \eqref{Chap2:eq2.2.3}
\begin{equation*}
\begin{split}
&\sup_{h\in\Omega_K}\int_{\|an\|<\delta_U}g^G_{t/2}((an)^{-1})|f(anh)-f(h)|dadn\\
&<\varepsilon\int_{\|an\|<\delta_U}g^G_{t/2}((an)^{-1})dadn\leq\varepsilon\int_A\int_N g^G_{t/2}((an)^{-1})dadn=\varepsilon,\\
&\sup_{h\in\Omega_K}\int_{\|an\|\geq\delta_U}g^G_{t/2}((an)^{-1})|f(anh)-f(h)|dadn\\
&\leq\sup_{h\in\Omega_K}\sup_{\|b\|\geq \delta_U}|f(bh)-f(h)|\int_{\|an\|\geq\delta_U}g^G_{t/2}((an)^{-1})dadn\\
&\leq2M\int_{\|an\|\geq\delta_U}g^G_{t/2}((an)^{-1})dadn<2M\varepsilon.
\end{split}
\end{equation*}
Thus we deduce from the inequality \eqref{SpecialFunctions:eq3.2.7} that
\begin{equation*}
\begin{split}
&\sup_{h\in\Omega_K}|g^{NA}_{t/2}*f^K_t(h)-f(h)|\\
&\leq\sup_{h\in\Omega_K}\int_A\int_Ng^G_{t/2}((an)^{-1})|f^K_t(anh)-f(anh)|dadn\\
&+\sup_{h\in\Omega_K}\int_{\|an\|<\delta_U}g^G_{t/2}((an)^{-1})|f(anh)-f(h)|dadn\\
&+\sup_{h\in\Omega_K}\int_{\|an\|\geq\delta_U}g^G_{t/2}((an)^{-1})|f(anh)-f(h)|dadn\to0
\end{split}
\end{equation*}
as $t$ approaches zero.

\noindent (iii) The maximum principle of the heat semigroup $e^{t\Delta_K}$ implies that $|f^K_t(h)|\leq M$ for $h\in G$. Hence for $h\in\Omega^c_K$, we have
\begin{equation*}
\begin{split}
&|g^{NA}_{t/2}*f^K_t(h)-f(h)|=|g^{NA}_{t/2}*f^K_t(h)|\leq\int_A\int_Ng^G_{t/2}((an)^{-1})|f^K_t(anh)|dadn\\
&\leq\int_A\int_Ng^G_{t/2}((an)^{-1})|f^K_t(anh)|dadn\leq M\int_A\int_Ng^G_{t/2}((an)^{-1})\mathbf{1}_{\Omega_K}(anh)dadn,
\end{split}
\end{equation*}
where $\mathbf{1}_{\Omega_K}$ denotes the indicator function of $\Omega_K$. Since $\mathbf{1}_{\Omega_K}\in L^1(G)\cap L^2(G)$ is bi-$K$-invariant, its spherical transform is given by $\hat{\mathbf{1}}_{\Omega_K}(\lambda)=\int_G\mathbf{1}_{\Omega_K}(x)\varphi_{\lambda}(x^{-1})dx$ for $\lambda\in i\mathfrak{a}^{*}$. Then we get the estimate
\begin{equation*}
\begin{split}
&M^{-2}\int_{\Omega^c_K}|g^{NA}_{t/2}*f^K_t(h)-f(h)|^2
dh\\
&\leq\int_{\Omega^c_K}\left(\int_A\int_Ng^G_{t/2}((an)^{-1})\mathbf{1}_{\Omega_K}(anh)dadn\right)^2dh\\
&=\int_{\Omega^c_K}\left|\int_A\int_Ng^G_{t/2}((an)^{-1})\mathbf{1}_{\Omega_K}(anh)dadn-\mathbf{1}_{\Omega_K}(h)\right|^2 dh\\
&=\int_{\Omega^c_K}\left|\int_K\int_A\int_Ng^G_{t/2}((kan)^{-1})\mathbf{1}_{\Omega_K}(kanh)dkdadn-\mathbf{1}_{\Omega_K}(h)\right|^2 dh\\
&=\int_{\Omega^c_K}\left|\int_Gg^G_{t/2}(x^{-1})\mathbf{1}_{\Omega_K}(xh)dx-\mathbf{1}_{\Omega_K}(h)\right|^2dh=\int_{\Omega^c_K}|g^G_{t/2}*\mathbf{1}_{\Omega_K}(h)-\mathbf{1}_{\Omega_K}(h)|^2dh\\
&\leq\int_G|g^G_{t/2}*\mathbf{1}_{\Omega_K}(h)-\mathbf{1}_{\Omega_K}(h)|^2dh=\frac{1}{|W|}\int_{i\mathfrak{a}^{*}}|\hat{\mathbf{1}}_{\Omega_K}(\lambda)|^2(e^{(\langle\lambda,\lambda\rangle-\langle\rho,\rho\rangle)\frac{t}{2}}-1)^2\frac{d\lambda}{|\mathbf{c}(\lambda)|^2}.
\end{split}
\end{equation*}
Here, in the last equality we use the Plancherel theorem and the fact that $g^G_t$ is the inverse spherical transform of $\lambda\mapsto e^{(\langle\lambda,\lambda\rangle-\langle\rho,\rho\rangle)t}$ on $i\mathfrak{a}^{*}$. We observe that $e^{(\langle\lambda,\lambda\rangle-\langle\rho,\rho\rangle)\frac{t}{2}}-1\to0$ pointwise on $i\mathfrak{a}^{*}$ as $t\to0+$. Moreover there exists a uniform bound
\begin{equation*}
|e^{(\langle\lambda,\lambda\rangle-\langle\rho,\rho\rangle)\frac{t}{2}}-1|\leq|e^{(\langle\lambda,\lambda\rangle-\langle\rho,\rho\rangle)\frac{t}{2}}|+1\leq2
\end{equation*}
for $\lambda\in i\mathfrak{a}^{*}$, and $\int_{i\mathfrak{a}^{*}}|\hat{\mathbf{1}}_{\Omega_K}(\lambda)|^2|\mathbf{c}(\lambda)|^{-2}d\lambda<\infty$. Hence by Lebesgue's dominated convergence theorem, the integral
\begin{equation*}
\int_{i\mathfrak{a}^{*}}|\hat{\mathbf{1}}_{\Omega_K}(\lambda)|^2(e^{(\langle\lambda,\lambda\rangle-\langle\rho,\rho\rangle)\frac{t}{2}}-1)^2\frac{d\lambda}{|\mathbf{c}(\lambda)|^2}
\end{equation*}
converges to $0$ as $t\to0+$. This proves that the integral \eqref{SpecialFunctions:eq3.2.6} converges to zero as $t\to0+$.
\end{proof}

\section{An integral expression of the heat kernel on a real group}\label{TheHeatKernelOnASemiSimpleGroup:sec4}
In this section we deduce Theorem~{\upshape\ref{TheHeatKernelOnASemiSimpleGroup:thm4.0.1}} concerning an integral expression of the heat kernel on $G_0$ whose integrand involves the heat kernel on $G$. We recall the Cartan involution $\sigma_0$ of the split real form $\mathfrak{g}_0$. Since $\mathfrak{g}_0$ is a split real form, $i\mathfrak{p}_0$ is a subspace of $\mathfrak{k}$. Let $\mathfrak{t}:=i\mathfrak{a}_0\subset i\mathfrak{p}_0$, then $\mathfrak{t}$ is maximal Abelian in $\mathfrak{k}$. Fix a Weyl chamber $\mathfrak{t}^{+}\subset \mathfrak{t}$. We denote by $\Sigma_K^{+}$ the set of positive restricted roots corresponding to $\mathfrak{t}^{+}$. We denote by $T$ the maximal torus subgroup of $K$ corresponding to $\mathfrak{t}$. Let $\overline{T^{+}}\subset T$ be the subset corresponding to the closure of $\mathfrak{t}^{+}$ via the exponential map. We have the Cartan decomposition $K=K_0\overline{T^{+}}K_0$.

Let
\begin{equation*}
\mathfrak{k}_\alpha=\{X\in \mathfrak{k}\oplus i\mathfrak{k}\ ;\ \operatorname{ad}(t)X= i\alpha(t)X \text{ for all } t\in\mathfrak{t}\}
\end{equation*}
denote the root space corresponding to $\alpha\in \Sigma_K^{+}$. For each restricted root $\alpha\in \Sigma_K^{+}$, let $m_\alpha=\dim_{\mathbb{R}}\mathfrak{k}_\alpha$ denote its multiplicity. Let
\begin{equation*}
\delta(\exp t):=
\prod_{\alpha\in\Sigma_K^{+}}\left|\sin(-i\alpha(t))\right|^{m_{\alpha}}
\end{equation*}
be the Jacobian determinant on $\mathfrak{t}$. We denote by $dk_0$ the normalized Haar measure on $K_0$ so that $\int_{K_0} dk_0=1$. Then the integration formula over the compact group $K$ using the decomposition $K=K_0\overline{T^{+}}K_0$ is given by
\begin{equation*}
\int_K f(k)dk=c_K\int_{K_0} \int_{\overline{T^{+}}} \int_{K_0} f(k_1\tau k^{-1}_1k_2) \delta(\tau)dk_1d\tau dk_2,
\end{equation*}
where $c_K$ is a normalization constant ensuring that the total volume integrates to $1$. We are ready to assert a key proposition on a relationship between $\rho_t^{K_0}$ and $\rho_t^K$.

\begin{prop}\label{TheHeatKernelOnASemiSimpleGroup:prop4.0.1}
The heat kernel $\rho_t^{K_0}$ on $K_0$ admits an integral expression written as
\begin{equation}\label{TheHeatKernelOnASemiSimpleGroup:eq4.0.1}
\rho_t^{K_0}(x)=(1-c_K)+c_K\int_{K_0}\int_{\overline{T^{+}}}\rho_{t/2}^K(k\tau k^{-1}x)\delta(\tau)dkd\tau.
\end{equation}
\end{prop}
\begin{proof}
Let $x\in K_0$, and let
\begin{equation*}
\nu^K_t(x):=\int_{K_0}\int_{\overline{T^{+}}}\rho_t^K(k\tau k^{-1}x)\delta(\tau)dkd\tau
\end{equation*}
denote a function on $K_0$. We first show that $\nu^K_{t/2}$ satisfies the heat equation on $K_0$. We see that
\begin{equation*}
\begin{split}
(\partial_t\nu^K_{t/2})(x)=\int_{K_0}\int_{\overline{T^{+}}}(\partial_t\rho_{t/2}^K)(k\tau k^{-1}x)\delta(\tau)dkd\tau=\frac{1}{2}\int_{K_0}\int_{\overline{T^{+}}}(\Delta_K\rho_{t/2}^K)(k\tau k^{-1}x)\delta(\tau)dkd\tau.
\end{split}
\end{equation*}
Using the decomposition $K\approx K/K_0\times K_0\approx (\operatorname{Ad}_K(K_0)\overline{T^{+}})K_0\times K_0$ induced by the Cartan decomposition of $\mathfrak{k}$, we obtain
\begin{equation*}
\begin{split}
\int_{K_0}\int_{\overline{T^{+}}}(\Delta_K\rho_{t/2}^K)(k\tau k^{-1}x)\delta(\tau)dkd\tau=2\int_{K_0}\int_{\overline{T^{+}}}(\Delta_{K_0}\rho_{t/2}^K)(k\tau k^{-1}x)\delta(\tau)dkd\tau.
\end{split}
\end{equation*}
Thus we deduce
\begin{equation*}
\begin{split}
(\partial_t\nu^K_{t/2})(x)=\int_{K_0}\int_{\overline{T^{+}}}(\Delta_{K_0}\rho_{t/2}^K)(k\tau k^{-1}x)\delta(\tau)dkd\tau,
\end{split}
\end{equation*}
which implies that $\nu^K_{t/2}$ is a solution to the heat equation on $K_0$. Hence $(1-c_K)+c_K\nu^K_{t/2}$ also solves the heat equation on $K_0$.

Recall that, since $\rho_{t/2}^K$ is the fundamental solution to the $L^2$-Cauchy problem, it is a smooth non-negative function on $(0,\infty)\times K$. Furthermore it follows from Lemma~{\upshape\ref{Chap1:lem1.1.2}} that, for any open neighborhood $V$ at $e$,
\begin{equation*}
\lim_{t\to0+}\int_V \rho_t^K(y)dy=1.
\end{equation*}
Hence, to see that $(1-c_K)+c_K\nu^K_{t/2}$ satisfies the initial condition, we show that, for any open neighborhood $U\subset K_0$ at $e$,
\begin{equation*}
\lim_{t\to0+}\int_U ((1-c_K)+c_K\nu^K_t(x))dx=1.
\end{equation*}
Take an open neighborhood $U\ni e$. Since
\begin{equation*}
(\operatorname{Ad}_K(K_0)\overline{T^{+}})U\approx(\operatorname{Ad}_K(K_0)\overline{T^{+}})K_0\times U\subset(\operatorname{Ad}_K(K_0)\overline{T^{+}})K_0\times K_0\approx K
\end{equation*}
is an open set containing $e$, we have
\begin{equation*}
\begin{split}
&\int_U ((1-c_K)+c_K\nu^K_t(x))dx=(1-c_K)+c_K\int_{K_0}\int_{\overline{T^{+}}}\int_U\rho_t^K(k\tau k^{-1}x)\delta(\tau)dkd\tau dx\\
&=(1-c_K)+c_K\int_{(\operatorname{Ad}_K(K_0)\overline{T^{+}})U}\rho_t^K(y)dy\to(1-c_K)+c_K=1
\end{split}
\end{equation*}
as $t\to0+$. Therefore by Lemma~{\upshape\ref{Chap1:lem1.1.2}} we obtain $(1-c_K)+c_K\nu^K_t\to\delta_e$ as $t\to0+$ in the distribution topology on $K_0$. This completes the proof.
\end{proof}

Using the reduction map $\mathbf{M}$, see the bijection \eqref{Chap1:eq1.3.4}, and Proposition~{\upshape\ref{TheHeatKernelOnASemiSimpleGroup:prop4.0.1}}, we deduce the main result of the present section.
\begin{theorem}\label{TheHeatKernelOnASemiSimpleGroup:thm4.0.1}
For $t>0$ and $g\in G_0$, the heat kernel $\rho^{G_0}_t(g)$ is expressed as
\begin{equation}\label{TheHeatKernelOnASemiSimpleGroup:eq4.0.2}
\begin{split}
\rho_t^{G_0}(g)=c_K\int_{K_0^{\mathbb{C}}}\int_{K_0}\int_{\overline{T^{+}}}\rho_{t/2}^G(ya_gk_gk\tau k^{-1})\delta(\tau)dkd\tau dy+(1-c_K)(\mathbf{M}g_{t/4}^G)(a_g),
\end{split}
\end{equation}
where $c_K$ is a constant determined by the Haar measure on $K$ and measures associated with the Cartan decomposition of $K$, and we put $a_g:=\operatorname{Crt}_{\overline{A^{+}}}(g)$ and $k_g:=\operatorname{Crt}_{K_0}(g)\operatorname{Crt}'_{K_0}(g)$.
\end{theorem}
\begin{proof}
For $\bullet\in\{0,\emptyset\}$, we recall the expression \eqref{SpecialFunctions:eq3.2.4} of the heat kernel $\rho_t^{G_{\bullet}}$ on $G_{\bullet}$
\begin{equation*}
\begin{split}
\rho_t^{G_{\bullet}}(h)&=[g_{t/2}^{G_{\bullet}}\times\rho_t^{K_{\bullet}}]([\operatorname{Iw}_{NA}(\sqrt{h(\sigma(h))^{-1}}),(\operatorname{Iw}_{K_{\bullet}}(\sqrt{h(\sigma(h))^{-1}}k_h))^{-1}])\\
&=\frac{1}{|W|}\int_{K_{\bullet}}\left(\int_{i\mathfrak{a}^{*}}e^{(\langle\lambda,\lambda\rangle-\langle\rho,\rho\rangle)\frac{t}{2}}e^{(\lambda+\rho)\log\operatorname{Iw}_A(x\sqrt{h(\sigma(h))^{-1}}x^{-1})}\frac{d\lambda}{|\mathbf{c}(\lambda)|^2}\right)\\
&\hspace{20mm}\times\rho^{K_{\bullet}}_t((\operatorname{Iw}_{K_{\bullet}}(x\sqrt{h(\sigma(h))^{-1}}x^{-1}k_h))^{-1})dx.
\end{split}
\end{equation*}
Then, for $h\in G_{\bullet}$, changing independent variables $x\mapsto x(\operatorname{Crt}_{K_{\bullet}}(h))^{-1}$ by right translation, we get the expression
\begin{equation*}
\rho_t^{G_{\bullet}}(h)=[g_{t/2}^{G_{\bullet}}\times\rho_t^{K_{\bullet}}]([\operatorname{Crt}_{\overline{A^{+}}}(h),(\operatorname{Crt}_{K_{\bullet}}(h)\operatorname{Crt}'_{K_{\bullet}}(h))^{-1}]).
\end{equation*}
Using the formulae \eqref{Chap1:eq1.3.3} and \eqref {Chap1:eq1.3.7}, we obtain
\begin{equation*}
\begin{split}
\rho_t^{G_0}(h)&=[g_{t/2}^{G_0}\times\rho_t^{K_0}]([\operatorname{Crt}_{\overline{A^{+}}}(h),(\operatorname{Crt}_{K_0}(h)\operatorname{Crt}'_{K_0}(h))^{-1}])\\
&=[\mathbf{M}g_{t/4}^G\times\rho_t^{K_0}]([\operatorname{Crt}_{\overline{A^{+}}}(h),(\operatorname{Crt}_{K_0}(h)\operatorname{Crt}'_{K_0}(h))^{-1}])\\
&=[\int_{K_0^{\mathbb{C}}} g_{t/4}^G(y\bullet)dy\times\rho_t^{K_0}]([\operatorname{Crt}_{\overline{A^{+}}}(h),(\operatorname{Crt}_{K_0}(h)\operatorname{Crt}'_{K_0}(h))^{-1}])\\
&=\int_{K_0^{\mathbb{C}}} g_{t/4}^G(\sqrt{ya_h(\sigma(ya_h))^{-1}})\rho_t^{K_0}((\operatorname{Iw}_{K_0}(\sqrt{ya_h(\sigma(ya_h))^{-1}})k_h)^{-1})dy.
\end{split}
\end{equation*}
By Proposition~{\upshape\ref{TheHeatKernelOnASemiSimpleGroup:prop4.0.1}} we derive
\begin{equation*}
\begin{split}
&[\int_{K_0^{\mathbb{C}}} g_{t/4}^G(y\bullet)dy\times\rho_t^{K_0}]([\operatorname{Crt}_{\overline{A^{+}}}(h),(\operatorname{Crt}_{K_0}(h)\operatorname{Crt}'_{K_0}(h))^{-1}])\\
&=c_K[\int_{K_0^{\mathbb{C}}} g_{t/4}^G(y\bullet)dy\times\int_{K_0}\int_{\overline{T^{+}}}\rho_{t/2}^K(k\tau k^{-1}\bullet)\delta(\tau)dkd\tau]([a_h,k_h^{-1}])\\
&+(1-c_K)[\int_{K_0^{\mathbb{C}}} g_{t/4}^G(y\bullet)dy\times\mathbf{1}]([a_h,k_h^{-1}])\\
&=c_K\int_{K_0^{\mathbb{C}}} g_{t/4}^G(\sqrt{ya_h(\sigma(ya_h))^{-1}})\\
&\times\int_{K_0}\int_{\overline{T^{+}}}\rho_{t/2}^K((k\tau k^{-1}\operatorname{Iw}_{K_0}(\sqrt{ya_h(\sigma(ya_h))^{-1}})k_h)^{-1})\delta(\tau)dkd\tau dy+(1-c_K)(\mathbf{M}g_{t/4}^G)(a_h),
\end{split}
\end{equation*}
where $\mathbf{1}$ denotes the constant function with value $1$ on $K$. We observe that
\begin{equation*}
\begin{split}
&\int_{K_0^{\mathbb{C}}} g_{t/4}^G(\sqrt{ya_h(\sigma(ya_h))^{-1}})\int_{K_0}\int_{\overline{T^{+}}}\rho_{t/2}^K((k\tau k^{-1}\operatorname{Iw}_{K_0}(\sqrt{ya_h(\sigma(ya_h))^{-1}})k_h)^{-1})\delta(\tau)dkd\tau dy\\
&=\int_{K_0^{\mathbb{C}}}\int_{K_0}\int_{\overline{T^{+}}} g_{t/4}^G(\sqrt{ya_h(\sigma(ya_h))^{-1}})\rho_{t/2}^K((\operatorname{Iw}_{K_0}(\sqrt{ya_h(\sigma(ya_h))^{-1}})k\tau k^{-1}k_h)^{-1})\delta(\tau)dkd\tau dy\\
&=[\int_{K_0^{\mathbb{C}}}g_{t/4}^G(y\bullet)dy\times\int_{K_0}\int_{\overline{T^{+}}}\rho_{t/2}^K(k\tau k^{-1}\bullet)\delta(\tau)dkd\tau](a_h,k_h^{-1})\\
&=\int_{K_0^{\mathbb{C}}}\int_{K_0}\int_{\overline{T^{+}}}\rho_{t/2}^G(ya_hk_hk\tau k^{-1})\delta(\tau)dkd\tau dy.
\end{split}
\end{equation*}
Therefore we deduce
\begin{equation*}
\rho_t^{G_0}(h)=c_K\int_{K_0^{\mathbb{C}}}\int_{K_0}\int_{\overline{T^{+}}}\rho_{t/2}^G(ya_hk_hk\tau k^{-1})\delta(\tau)dkd\tau dy+(1-c_K)(\mathbf{M}g_{t/4}^G)(a_h).
\end{equation*}
This proves Theorem~{\upshape\ref{TheHeatKernelOnASemiSimpleGroup:thm4.0.1}}.
\end{proof}

We finally point out that the relationship \eqref{TheHeatKernelOnASemiSimpleGroup:eq4.0.1} between $\rho^{G_0}_t$ and $\rho^{ G}_t$ arises through the push-forward of $\rho^{ G}_t$ under projection along fibres. Let
\begin{equation*}
p_{ G}:\ G\to G_0,\ nax_1\tau x_2\mapsto ax_1x_2
\end{equation*}
be a projection associated with the Iwasawa decomposition of $ G$ and the Cartan decomposition of $ K$, where $n\in N$, $\tau\in \overline{T^{+}}$, $a\in A$ and $x_1,x_2\in K_0$. Recall that $\tilde{\sigma}_0$ and the Cartan involution $\sigma$ commute. Then one deduces the product decomposition of the symmetric space
\begin{equation}\label{Chap2:eq3.0.3}
G/ K\approx G^{\tilde{\sigma}_0}/K_0\times G_0/K_0,
\end{equation}
see \cite[Section 3]{Ross79}. Thus the map $\mathbf{M}$ in the identity \eqref{SpecialFunctions:eq4.3.3} realizes integration over spaces traced by $G^{\tilde{\sigma}_0}$-orbits of points in $G_0/K_0$. Recall that, to deduce the integral expression \eqref{SpecialFunctions:eq4.3.2} for the heat kernel $\rho^{G_0}_t$, we employ Proposition~{\upshape\ref{TheHeatKernelOnASemiSimpleGroup:prop4.0.1}} concerning an integral formula for $\rho^{K_0}_t$. For our purpose, the integral \eqref{TheHeatKernelOnASemiSimpleGroup:eq4.0.2} of the heat kernel $\rho^{ G}_t$ over $\operatorname{Ad}_K(K_0)\overline{T^{+}}$ is regarded as integration over spaces traced by $ K/K_0$-fibres of points in $K_0$. Thus we interpret the formula \eqref{TheHeatKernelOnASemiSimpleGroup:eq4.0.2} for $\rho^{G_0}_t$ as the push-forward of $\rho^{ G}_t$ under $p_{ G}$, namely, as an integral over the fibers.

\section{An example: the case of $SL(2,\mathbb{R})$}\label{TheHeatKernelOnASemiSimpleGroup:sec5}
Fix $t>0$. In this section we assume $ G=SL(2,\mathbb{C})$ $ K=SU(2)$, $G_0=SL(2,\mathbb{R})$ and $K_0=SO(2)$.

\subsection{Main results for $SL(2,\mathbb{R})$}
Let $\sigma_0$ be the Cartan involution on $G_0$ given by $\sigma_0(g):=({}^tg)^{-1}$ for $g\in G_0$. Let $\mathfrak{g}_0$ and $\mathfrak{k}_0$ denote the Lie algebras $\mathfrak{sl}(2,\mathbb{R})$ and $\mathfrak{so}(2)$ respectively. We take a basis of $\mathfrak{g}_0$
\begin{equation*}
\left\{Z_1:=\frac{1}{2}\begin{pmatrix}
1 & 0\\
0 & -1
\end{pmatrix},\ 
Z_2:=\frac{1}{2}\begin{pmatrix}
0 & 1\\
1 & 0
\end{pmatrix},\ Z_3:=
\frac{1}{2}\begin{pmatrix}
0 & -1\\
1 & 0
\end{pmatrix}\right\}.
\end{equation*}
Let $\mathfrak{p}_0:=\operatorname{Span}_{\mathbb{R}}\{Z_1,Z_2\}$. We have $\mathfrak{k}_0=\mathbb{R}Z_3$ and the Cartan decomposition $\mathfrak{g}_0=\mathfrak{p}_0\oplus\mathfrak{k}_0$ for $\sigma_0$. We define a Riemannian metric $(g_{ij})$ on $G_0$ by
\begin{equation*}
(g_{ij})_h:=g_h((L_h)_{*}Z_i,(L_h)_{*}Z_j)=\delta_{ij}
\end{equation*}
with $h\in G_0$ and $1\leq i,j\leq3$, where $L_h$ is the left translation $x\mapsto hx$ of $G_0$. Given $Z\in\mathfrak{g}_0$, we denote the left-$G_0$-invariant vector field corresponding to $Z$ by
\begin{equation*}
(\tilde{Z}f)(g):=\left.\DS\frac{d}{ds}\right|_{s=0}f(g(\exp sZ)),
\end{equation*}
where $g\in G_0$ and $f$ denotes a compactly supported smooth function on $G_0$. Recall that we can uniquely extend $\tilde{Z}$ to an operator on $L^2(G_0,dg)$. We denote by $\Delta_{G_0}$ the Laplace--Beltrami operator associated with $(g_{ij})$, which are written as
\begin{equation*}
-2\Delta_{G_0}=\tilde{Z}_1\circ\tilde{Z}_1+\tilde{Z}_2\circ\tilde{Z}_2+\tilde{Z}_3\circ\tilde{Z}_3.
\end{equation*}
By Definition~{\upshape\ref{SpecialFunctions:def2.1.1}}, we obtain the heat kernel $\rho^{G_0}_t$ with respect to the Riemannian metric $(g_{ij})$.\\

We now deduce an expression for the heat kernel $\rho^{G_0}_t$ in relation to the heat Gaussian $g^{G_0}_t$ and the heat kernel $\rho^{K_0}_t$, using the argument from the proof of Theorem~{\upshape\ref{SpecialFunctions:thm3.0.1}}. In the case of $G_0=SL(2,\mathbb{R})$, it is known that $|\mathbf{c}(i\nu)|^{-2}=\pi\nu\tanh(\pi\nu)$ for $\nu\in\mathbb{R}$. Hence, by Proposition~{\upshape\ref{SpecialFunctions:prop2.1.1}}, we get an explicit expression of $g^{G_0}_t$
\begin{equation*}
g^{G_0}_t(x)=\frac{1}{2}\int_{\mathbb{R}}e^{-\left(\nu^2+\frac{1}{4}\right)t}\varphi_{i\nu}(x)\frac{\nu\tanh(\pi \nu)d\nu}{2\pi}
\end{equation*}
for $x\in G_0$. Let $(S^1,\theta)$ be the unit circle with the Laplace operator $\tfrac{1}{2}\tfrac{\partial^2}{\partial\theta^2}$. As $K_0=\{k_{\theta}:=\left(\begin{smallmatrix}
\cos\theta & -\sin\theta\\
\sin\theta & \cos\theta
\end{smallmatrix}\right)\ ;\ \theta\in\mathbb{R}\}\simeq S^1$, we have
\begin{equation*}
\rho^{K_0}_t(k_{\theta})=\frac{1}{2\pi}\vartheta\left(\frac{\theta}{2\pi},\frac{t}{2\pi}i\right),
\end{equation*}
where $\vartheta(z,\tau):=\sum_{m\in\mathbb{Z}}e^{m^2i\pi\tau}e^{2m\pi iz}$ is Jacobi's theta function for $z,\tau\in\mathbb{C}$ with $\Im\tau>0$. In this setup, specializing Theorem~{\upshape\ref{SpecialFunctions:thm3.0.1}} yields an explicit form \eqref{SpecialFunctions:eq4.0.1} of the heat kernel $\rho^{G_0}_t$ in Corollary~{\upshape\ref{SpecialFunctions:cor4.0.1}}.

\begin{cor}\label{SpecialFunctions:cor4.0.1}
Let $\{\alpha^0\}$ be a positive root system of $\mathfrak{g}_0=\mathfrak{sl}(2,\mathbb{R})$. Then for $t>0$ and $g\in G_0$, the heat kernel $\rho^{G_0}_t$ is expressed as
\begin{equation}\label{SpecialFunctions:eq4.0.1}
\begin{split}
\rho^{G_0}_t(g)=\frac{1}{2}\int_{K_0}\int_{\mathbb{R}}&e^{-\left(\nu^2+\frac{1}{4}\right)\frac{t}{2}}e^{\left(i\nu+\frac{1}{2}\right)\alpha^0(\log\operatorname{Iw}_A(k\sqrt{g(\sigma(g))^{-1}}k^{-1}))}\frac{\nu\tanh(\pi \nu)d\nu}{2\pi}\\
&\times\rho^{K_0}_t((\operatorname{Iw}_{K_0}(k\sqrt{g(\sigma(g))^{-1}}k^{-1})\operatorname{Crt}_{K_0}(g)\operatorname{Crt}'_{K_0}(g))^{-1})dk,
\end{split}
\end{equation}
where we put $\sqrt{g(\sigma(g))^{-1}}:=\operatorname{Crt}_{K_0}(g)\operatorname{Crt}_{\overline{A^{+}}}(g)(\operatorname{Crt}_{K_0}(g))^{-1}$.
\end{cor}

\subsection{An integral presentation of $\rho_t^{SL(2,\mathbb{R})}$ via Tchebycheff polynomials}
As a special case of Theorem~{\upshape\ref{TheHeatKernelOnASemiSimpleGroup:thm4.0.1}}, we deduce Proposition~{\upshape\ref{SpecialFunctions:prop4.3.1}} regarding an integral presentation of $\rho_t^{G_0}$ by means of Tchebycheff polynomials. We shall derive the expression \eqref{SpecialFunctions:eq4.3.1} of the heat kernel $\rho^{ K}_t$, using the Tchebycheff polynomials of the second kind. Let $T:=\{t_{\varphi}:=\operatorname{diag}[e^{i\varphi},e^{-i\varphi}]\ ;\ \varphi\in\mathbb{R}\}\subset K$ be a maximal torus, and $\Lambda_{\mathrm{dom}}:=\{l\ ;\ 2l\in\mathbb{Z},l\geq0\}$ be the set of dominant weights for the pair $( K,T)$. For each $l\in\Lambda_{\mathrm{dom}}$, there is a character $\chi_l$ given by the Weyl character formula
\begin{equation*}
\tau_l(t_{\varphi})=\frac{\sin (2l+1)\varphi}{\sin\varphi}=U_{2l}(\cos\varphi),
\end{equation*}
where $U_{2l}$ is the $(2l)$-th Tchebycheff polynomial of the second kind. $\chi_l$ is the trace of the irreducible unitary representation $\tau_l$ on the representation space $V_l$ with $\dim V_l=2l+1$. Then the eigenvalues of the Laplace--Beltrami operator $\Delta_{ K}$ are given by
\begin{equation*}
\Delta_{ K}\chi_l=-\lambda_l\chi_l,\ \lambda_l=\frac{l(l+1)}{2}.
\end{equation*}
For $k\in K$, since the $T$-component of $k$ for the Cartan decomposition $K=K_0TK_0$ is determined up to permutation and can be written as
\begin{equation*}
t_{\operatorname{Arc}\!\cos(\operatorname{tr}(k)/2)}=\operatorname{diag}[e^{i\operatorname{Arc}\!\cos(\operatorname{tr}(k)/2)},e^{-i\operatorname{Arc}\!\cos(\operatorname{tr}(k)/2)}]\in T,
\end{equation*} 
Theorem~{\upshape\ref{SpecialFunctions:thm2.1.2}} implies that we have
\begin{equation*}
\rho^{ K}_t(k)=\sum_{l\in\Lambda_{\mathrm{dom}}}(2l+1)e^{-\frac{l(l+1)}{2}t}U_{2l}(\operatorname{tr}(k)/2).
\end{equation*}
Here, given $g\in G$, $\operatorname{tr}(g)$ is the trace of $g$. In particular, for $k=\left(\begin{smallmatrix}
\alpha & \beta\\
-\bar{\beta} & \bar{\alpha}
\end{smallmatrix}\right)
\in K=SU(2)$, we have $\operatorname{tr}(k)=\alpha+\bar{\alpha}$. Setting $m=2l+1$ yields
\begin{equation}\label{SpecialFunctions:eq4.3.1}
\rho^{ K}_t(k)=\sum^{\infty}_{m=1}me^{-(m^2-1)\frac{t}{8}}U_{m-1}(\operatorname{tr}(k)/2).
\end{equation}

It is well-known that the $m$-th Tchebycheff polynomial $T_m$ of the first kind satisfies the derivative formula $T'_m=mU_{m-1}$ for integers $m\geq1$. Hence the heat kernel $\rho^{K_0}_t$ can be written as
\begin{equation*}
\begin{split}
&\rho^{K_0}_t(k_{\theta})=\frac{1}{2\pi}\vartheta\left(\frac{\theta}{2\pi},\frac{t}{2\pi}i\right)=\frac{1}{2\pi}+\frac{1}{\pi}\sum^{\infty}_{m=1}e^{-\frac{m^2}{2}t}T_m(\operatorname{tr}(k_{\theta})/2)\\
&=\frac{1}{2\pi}+\frac{1}{\pi}\sum^{\infty}_{m=1}e^{-\frac{m^2}{2}t}\int^{\frac{\pi}{2}}_0 mU_{m-1}(\operatorname{tr}(k_{\theta}\operatorname{diag}[e^{i\varphi},e^{-i\varphi}])/2)\\
&\times\sqrt{(\operatorname{tr}(k_{\theta})/2)^2-(\operatorname{tr}(k_{\theta}\operatorname{diag}[e^{i\varphi},e^{-i\varphi}]))/2)^2}d\varphi\\
&=\frac{1}{2\pi}+\frac{1}{\pi}\sum^{\infty}_{m=1}e^{-\frac{m^2}{2}t}\int^{\frac{\pi}{2}}_0 mU_{m-1}(\operatorname{tr}(\operatorname{diag}[e^{-i\varphi},e^{i\varphi}]k_{\theta})/2)\\
&\times\sqrt{(\operatorname{tr}(k_{\theta})/2)^2-(\operatorname{tr}(k_{\theta}\operatorname{diag}[e^{i\varphi},e^{-i\varphi}])/2)^2}d\varphi
\end{split}
\end{equation*}
for $k_{\theta}:=\left(\begin{smallmatrix}
\cos\theta & -\sin\theta\\
\sin\theta & \cos\theta
\end{smallmatrix}\right)\in K_0$ with $\theta\in\mathbb{R}$. Recall that the bilinear forms $B_0$ and $B$ respectively induces the metric on $K_0$ and that on $K$. Using the expression \eqref{SpecialFunctions:eq4.3.1} of $\rho^{ K}_t$, we get
\begin{equation*}
\begin{split}
\rho^{K_0}_t(k_{\theta})&=\frac{1}{2\pi}+\frac{e^{\frac{t}{4}}}{\pi}\int^{\frac{\pi}{2}}_0\rho^{ K}_{2t}(\operatorname{diag}[e^{-i\varphi},e^{i\varphi}]k_{\theta})\\
&\times\sqrt{(\operatorname{tr}(k_{\theta})/2)^2-(\operatorname{tr}(k_{\theta}\operatorname{diag}[e^{i\varphi},e^{-i\varphi}])/2)^2}d\varphi,
\end{split}
\end{equation*}
where the time variable $t$ in the heat kernel on $K$ is multiplied by $1/2$, which arises from the equation $B_
{\mathfrak{k}}=2(B_0)_
{\mathfrak{k}_0}$. Hence taking the coordinate $\theta/2$ of $K$ for $\theta\in[0,4\pi)$ yields
\begin{equation}\label{SpecialFunctions:eq4.3.2}
\begin{split}
\rho^{K_0}_t(k_{\theta/2})
&=\frac{1}{4\pi}+\frac{e^{\frac{t}{4}}}{2\pi}\int^{\frac{\pi}{2}}_0\rho^{ K}_{t/2}(\operatorname{diag}[e^{-i\varphi},e^{i\varphi}]k_{\theta/2})\\
&\times\sqrt{(\operatorname{tr}(k_{\theta/2})/2)^2-(\operatorname{tr}(k_{\theta/2}\operatorname{diag}[e^{i\varphi},e^{-i\varphi}])/2)^2}d\varphi.
\end{split}
\end{equation}
Note that, in the case of $K=SU(2)$, we have $\rho_t^{SU(2)}(k\tau k^{-1}k_{\theta/2})=\rho_t^{SU(2)}(\tau k^{-1}k_{\theta/2}k)=\rho_t^{SU(2)}(\tau k_{\theta/2})$ for $\tau\in T$ and $k,k_{\theta/2}\in K_0=SO(2)$. Cf. Proposition~{\upshape\ref{TheHeatKernelOnASemiSimpleGroup:prop4.0.1}}.

Next, we describe the heat Gaussian $g^{G_0}_t$ on $G_0$ by using $g^{ G}_t$ on $ G$. We follow part of the arguments in \cite[Section 10]{Fl-J78}. After a change of variables, $g^{ G}_t$ and $g^{G_0}_t$ take the form
\begin{equation*}
\begin{split}
g^{ G}_t(a_{r/2})=\frac{e^{-t}}{(4\pi t)^{\frac{3}{2}}}\frac{re^{-\frac{r^2}{4t}}}{\sinh r},\ g^{G_0}_t(a_{r/2})=\frac{\sqrt{2}e^{-\frac{t}{4}}}{(4\pi t)^{\frac{3}{2}}}\int^{\infty}_r\frac{se^{-\frac{s^2}{4t}}ds}{\sqrt{\cosh s-\cosh r}}.
\end{split}
\end{equation*}
Cf. \cite[Theorem 2.5.5]{JL08} for the scaling property of $g^G_t$ with respect to the dependence on the time variable $t$, and the dependence of the coefficient on the bilinear form $B$. For spherical functions $\Phi\in \mathscr{S}^2( K\backslash G/ K)$, the bijective linear operator \eqref{Chap1:eq1.3.4} of $\mathscr{S}^2( K\backslash G/ K)$ onto $\mathscr{S}^2(K_0\backslash G_0/K_0)$ for $n=2$, denoted by $\mathbf{M}$, can be explicitly written as
\begin{equation}\label{SpecialFunctions:eq4.3.5}
(\mathbf{M}\Phi)(a_r)=2^{-1}\int^{\infty}_{s=r}\frac{\Phi(a_{s/2})}{\sqrt{\cosh(2s)-\cosh(2r)}}\frac{d(\cosh^2s)}{\cosh s},
\end{equation}
which was initiated by Flensted-Jensen \cite{Fl-J78}. Cf. \cite[Equation 10.3]{Fl-J78}, where it is defined for the Killing form on $\mathfrak{g}_0$ and that on $\mathfrak{g}$, while we take bilinear forms $B_0$ and $B$ defined as $2$ times the Killing forms respectively. Then we obtain
\begin{equation}\label{SpecialFunctions:eq4.3.3}
g^{G_0}_{t/2}(a_{\frac{r}{2}})=(\mathbf{M}g^{ G}_{t/4})(a_{\frac{r}{2}}).
\end{equation}
Cf. \cite[Equations (6.20) and (10.6)]{Fl-J78}.\\

We come to the main result of this subsection. With respect to the Cartan decomposition $K=K_0TK_0$, we select a fundamental domain
\begin{equation*}
T_{\mathrm{fund}}:=\{t_{\varphi}=\operatorname{diag}[e^{i\varphi},e^{-i\varphi}]\ ;\ 0\leq\varphi<\pi/2\}\subset T
\end{equation*}
for the action of $\{\pm e,k_{\pm\pi/2}\}\subset K_0$. Let $d\tau$ denote the restriction of the Haar measure on $T\subset K$ to $T_{\mathrm{fund}}$.
\begin{prop}\label{SpecialFunctions:prop4.3.1}
For $t>0$ and $g\in G$, the heat kernel $\rho^G_t(g)$ is expressed as
\begin{equation}\label{SpecialFunctions:eq4.3.4}
\begin{split}
&\rho^{G_0}_t(g)\\
&=2^{-\frac{5}{2}}\pi^{-1}e^{\frac{t}{4}}\int^{\infty}_{s=\|\operatorname{Crt}_{\overline{A^{+}}}(g)\|}\int_{T_{\mathrm{fund}}}\sqrt{\frac{(\operatorname{tr}(\operatorname{Crt}_K(g)\operatorname{Crt}'_K(g)))^2-(\operatorname{tr}(\operatorname{Crt}_K(g)\operatorname{Crt}'_K(g)\tau))^2}{(\operatorname{tr}(a_s))^2-(\operatorname{tr}(\operatorname{Crt}_{\overline{A^{+}}}(g)))^2}}\\
&\rho^{ G}_{t/2}(a_{s/2}\operatorname{Crt}_K(g)\operatorname{Crt}'_K(g)\tau)d\tau\frac{d(\cosh^2s)}{\cosh s}+\frac{\mathbf{M}g^{ G}_{t/4}(\operatorname{Crt}_{\overline{A^{+}}}(g))}{4\pi}.
\end{split}
\end{equation}
\end{prop}
\begin{proof}
Let $g=k_{\theta_1/2}a_{r/2}k_{\theta_2/2}\in G_0=K_0\overline{A^{+}}K_0$. By Theorem~{\upshape\ref{SpecialFunctions:thm3.0.1}} we get
\begin{equation*}
\begin{split}
\rho^{G_0}_t(g)=\frac{1}{2}\int^{4\pi}_0\int_{\mathbb{R}}&e^{-\left(\nu^2+\frac{1}{4}\right)\frac{t}{2}}e^{\left(\nu i-\frac{1}{2}\right)\alpha^0(\log\operatorname{Iw}_A(k_{\frac{\theta+\theta_1}{2}}a_{r/2}k^{-1}_{\frac{\theta+\theta_1}{2}}))}\frac{\nu\tanh(\pi \nu)d\nu}{2\pi}\\
&\times\rho^{K_0}_t((\operatorname{Iw}_{K_0}(k_{\frac{\theta+\theta_1}{2}}a_{r/2}k^{-1}_{\frac{\theta+\theta_1}{2}})k_{\frac{\theta_1+\theta_2}{2}})^{-1})\frac{d\theta}{4\pi}\\
=\frac{1}{2}\int^{4\pi}_0\int_{\mathbb{R}}&e^{-\left(\nu^2+\frac{1}{4}\right)\frac{t}{2}}e^{\left(\nu i-\frac{1}{2}\right)\alpha^0(\log\operatorname{Iw}_A(k_{\theta/2}a_{r/2}k_{-\theta/2}))}\frac{\nu\tanh(\pi \nu)d\nu}{2\pi}\\
&\times\rho^{K_0}_t((\operatorname{Iw}_{K_0}(k_{\theta/2}a_{r/2}k_{-\theta/2})k_{\frac{\theta_1+\theta_2}{2}})^{-1})\frac{d\theta}{4\pi},
\end{split}
\end{equation*}
and, for $s\geq0$ and $\tau\in T$,
\begin{equation*}
\begin{split}
\rho^G_t(k_{\theta_1/2}a_{s/2}&k_{\theta_2/2}\tau)\\
=\frac{1}{2}\int_K\int_{\mathbb{R}}&e^{-\left(\nu^2+1\right)\frac{t}{2}}e^{\left(\nu i-1\right)\alpha^0(\log\operatorname{Iw}_A(kk_{\theta_1/2}a_{s/2}(kk_{\theta_1/2})^{-1}))}\frac{4\sinh^2\nu}{\langle\rho,\alpha^0\rangle^2}\frac{\sqrt{2}d\nu}{\pi}\\
&\times\rho^K_t((\operatorname{Iw}_K(kk_{\theta_1/2}a_{s/2}(kk_{\theta_1/2})^{-1})k_{\frac{\theta_1+\theta_2}{2}}\tau)^{-1})dk\\
=\frac{1}{2}\int_K\int_{\mathbb{R}}&e^{-\left(\nu^2+1\right)\frac{t}{2}}e^{\left(\nu i-1\right)\alpha^0(\log\operatorname{Iw}_A(ka_{s/2}k^{-1}))}\frac{4\sqrt{2}\sinh^2\nu d\nu}{\pi}\\
&\times\rho^K_t(\tau^{-1}k_{-\frac{\theta_1+\theta_2}{2}}(\operatorname{Iw}_K(ka_{s/2}k^{-1}))^{-1})dk.
\end{split}
\end{equation*}
Here $\alpha^0$ denotes the positive restricted root in $\Sigma^{+}=\Sigma^{+}_0$. Recall that the bilinear forms $B_0$ and $B$ on $\mathfrak{a}_0=\mathfrak{a}$ induce the polar heights $\|\cdot\|_0$ and $\|\cdot\|$ on $A$ respectively. Since $\|a_r\|^2_0=2\|a_r\|^2$ for $r\geq0$, the formula \eqref{SpecialFunctions:eq4.3.3} for $g^{G_0}_t$ and $g^G_t$ takes the explicit form
\begin{equation*}
\begin{split}
&\frac{1}{2}\int^{4\pi}_0\int_{\mathbb{R}}e^{-\left(\nu^2+\frac{1}{4}\right)\frac{t}{2}}e^{\left(\nu i-\frac{1}{2}\right)\alpha^0(\log\operatorname{Iw}_A(k_{\theta/2}a_{r/2}k_{-\theta/2}))}\frac{\nu\tanh(\pi \nu)d\nu}{2\pi}\frac{d\theta}{4\pi}\\
&=2^{-1}\int^{\infty}_{s=\frac{r}{2}}\frac{g^G_{\frac{t}{4}}(a_{s/2})}{\sqrt{\cosh(2s)-\cosh r}}\frac{d(\cosh^2s)}{\cosh s}.
\end{split}
\end{equation*}
Hence, via the diffeomorphism $G\approx P$, we obtain the identity
\begin{equation*}
[g^{G_0}_{t/2}\times\rho^{K_0}_t]([a_{r/2},k_{-\frac{\theta_1+\theta_2}{2}}])=[\mathbf{M}g^G_{t/4}\times\rho^{K_0}_t]([a_{r/2},k_{-\frac{\theta_1+\theta_2}{2}}]),
\end{equation*}
which can be written as 
\begin{equation}\label{Chap2:eq3.0.4}
\begin{split}
&\frac{1}{2}\int^{4\pi}_0\int_{\mathbb{R}}e^{-\left(\nu^2+\frac{1}{4}\right)\frac{t}{2}}e^{\left(\nu i-\frac{1}{2}\right)\alpha^0(\log\operatorname{Iw}_A(k_{\theta/2}a_{r/2}k_{-\theta/2}))}\frac{\nu\tanh(\pi \nu)d\nu}{2\pi}\\
&\times\rho^{K_0}_t((\operatorname{Iw}_{K_0}(k_{\theta/2}a_{r/2}k_{-\theta/2})k_{\frac{\theta_1+\theta_2}{2}})^{-1})\frac{d\theta}{4\pi}\\
&=\frac{1}{2}\int^{\infty}_{s=\frac{r}{2}}\frac{[g^G_{\frac{t}{4}}\times\rho^{K_0}_t]([a_{s/2},k_{-\frac{\theta_1+\theta_2}{2}}])}{\sqrt{\cosh(2s)-\cosh r}}\frac{d(\cosh^2s)}{\cosh s}.
\end{split}
\end{equation}

We shall deduce the formula \eqref{SpecialFunctions:eq4.3.4} for $g=k_{\theta_1/2}a_{r/2}k_{\theta_2/2}$. Using the equation \eqref{Chap2:eq3.0.4}, we compute
\begin{equation*}
\begin{split}
&\rho^{G_0}_t(g)=\rho^{G_0}_t(k_{\theta_1/2}a_{r/2}k_{\theta_2/2})\\
&=\frac{1}{2}\int^{4\pi}_0\left(\int_{\mathbb{R}}e^{-\left(\nu^2+\frac{1}{4}\right)\frac{t}{2}}e^{\left(\nu i-\frac{1}{2}\right)\alpha^0(\log\operatorname{Iw}_A(k_{\theta/2}a_{r/2}k_{-\theta/2}))}\frac{\nu\tanh(\pi \nu)d\nu}{2\pi}\right)\\
&\times\rho^{K_0}_t((\operatorname{Iw}_{K_0}(k_{\theta/2}a_{r/2}k_{-\theta/2})k_{\frac{\theta_1+\theta_2}{2}})^{-1})\frac{d\theta}{4\pi}\\
&=\frac{1}{2}\int^{\infty}_{s=\frac{r}{2}}\frac{[g^G_{\frac{t}{4}}\times\rho^{K_0}_t]([a_{s/2},k_{-\frac{\theta_1+\theta_2}{2}}])}{\sqrt{\cosh(2s)-\cosh r}}\frac{d(\cosh^2s)}{\cosh s}.
\end{split}
\end{equation*}
Combining with the equation \eqref{SpecialFunctions:eq4.3.2}, we get
\begin{equation*}
\begin{split}
&[g^G_{\frac{t}{4}}\times\rho^{K_0}_t]([a_{s/2},k_{-\frac{\theta_1+\theta_2}{2}}])\\
&=\frac{e^{\frac{t}{4}}}{2\pi}\int_{T_{\mathrm{fund}}}[g^G_{\frac{t}{4}}\times \rho^K_{t/2}]([a_{s/2},\tau^{-1} k_{-\frac{\theta_1+\theta_2}{2}}])\sqrt{(\operatorname{tr}(k_{-\frac{\theta_1+\theta_2}{2}})/2)^2-(\operatorname{tr}(k_{-\frac{\theta_1+\theta_2}{2}}\tau)/2)^2}d\tau\\
&+\frac{1}{4\pi}[g^G_{\frac{t}{4}}\times \mathbf{1}]([a_{s/2},k_{-\frac{\theta_1+\theta_2}{2}}]).
\end{split}
\end{equation*}
for $s\geq r/2$, where $\mathbf{1}$ denotes the constant function $k\mapsto 1$ on $K$. Note that
\begin{equation*}
\begin{split}
&[g^G_{\frac{t}{4}}\times \mathbf{1}]([a_{s/2},k_{-\frac{\theta_1+\theta_2}{2}}])\\
&=\frac{1}{2}\int_K\int_{\mathbb{R}}e^{-\left(\nu^2+1\right)\frac{t}{4}}e^{\left(\nu i-1\right)\alpha^0(\log\operatorname{Iw}_A(ka_{s/2}k^{-1}))}\frac{4\sqrt{2}\sinh^2\nu d\nu}{\pi}\\
&\times\mathbf{1}((\operatorname{Iw}_K(ka_{s/2}k^{-1})k_{\frac{\theta_1+\theta_2}{2}})^{-1})dk\\
&=\frac{1}{2}\int_{\mathbb{R}}e^{-\left(\nu^2+1\right)\frac{t}{4}}\left(\int_Ke^{\left(\nu i-1\right)\alpha^0(\log\operatorname{Iw}_A(ka_{s/2}k^{-1}))}dk\right)\frac{4\sqrt{2}\sinh^2\nu d\nu}{\pi}=g^G_{\frac{t}{4}}(a_{s/2}).
\end{split}
\end{equation*}
Since we observe that
\begin{equation*}
\begin{split}
&[g^G_{\frac{t}{4}}\times \rho^K_{t/2}]([a_{s/2},\tau^{-1}k_{-\frac{\theta_1+\theta_2}{2}}])\\
&=\frac{1}{2}\int_K\int_{\mathbb{R}}e^{-\left(\nu^2+1\right)\frac{t}{4}}e^{\left(\nu i-1\right)\alpha^0(\log\operatorname{Iw}_A(ka_{s/2}k^{-1}))}\frac{4\sqrt{2}\sinh^2\nu d\nu}{\pi}\\
&\times\rho^K_t(\tau^{-1}k_{-\frac{\theta_1+\theta_2}{2}}(\operatorname{Iw}_K(ka_{s/2}k^{-1}))^{-1})dk=\rho^G_{t/2}(a_{s/2}k_{\frac{\theta_1+\theta_2}{2}}\tau)
\end{split}
\end{equation*}
for $s\geq r/2$ and $\tau\in T_{\mathrm{fund}}$, we deduce the integral formula \eqref{SpecialFunctions:eq4.3.4} for $\rho^{G_0}_t$
\begin{equation*}
\begin{split}
&\rho^{G_0}_t(g)=\rho^{G_0}_t(k_{\theta_1/2}a_{r/2}k_{\theta_2/2})\\
&=\frac{1}{2}\frac{e^{\frac{t}{4}}}{2\pi}\int^{\infty}_{s=\frac{r}{2}}\int_{T_{\mathrm{fund}}}\rho^G_{t/2}(a_{s/2}k_{\frac{\theta_1+\theta_2}{2}}\tau)\sqrt{\frac{(\operatorname{tr}(k_{-\frac{\theta_1+\theta_2}{2}})/2)^2-(\operatorname{tr}(k_{-\frac{\theta_1+\theta_2}{2}}\tau)/2)^2}{\cosh(2s)-\cosh r}}d\tau\frac{d(\cosh^2s)}{\cosh s}\\
&+\frac{1}{4\pi}\left(\frac{1}{2}\int^{\infty}_{s=\frac{r}{2}}\frac{g^G_{t/4}(a_{s/2})}{\sqrt{\cosh(2s)-\cosh r}}\frac{d(\cosh^2s)}{\cosh s}\right)\\
&=2^{-\frac{3}{2}}\frac{e^{\frac{t}{4}}}{2\pi}\int^{\infty}_{s=\frac{r}{2}}\int_{T_{\mathrm{fund}}}\rho^G_{t/2}(a_{s/2}k_{\frac{\theta_1+\theta_2}{2}}\tau)\sqrt{\frac{(\operatorname{tr}(k_{\frac{\theta_1+\theta_2}{2}})/2)^2-(\operatorname{tr}(k_{\frac{\theta_1+\theta_2}{2}}\tau)/2)^2}{\cosh^2s-\cosh^2 (r/2)}}d\tau\frac{d(\cosh^2s)}{\cosh s}\\
&+\frac{\mathbf{M}g^G_{t/4}(a_{r/2})}{4\pi}
\end{split}
\end{equation*}
if $g=k_{\theta_1/2}a_{r/2}k_{\theta_2/2}\in G_0$.
\end{proof}

\subsection{A formula for the heat kernel $\rho_t^{SL(2,\mathbb{R})}$ in previous works}
A previous approach to computing $\rho_t^{G_0}$ uses the fibration $K_0\to G_0\to G_0/K_0$ via the heat kernel on the universal covering group $\widetilde{SL(2,\mathbb{R})}$ of $G_0=SL(2,\mathbb{R})$, as presented in \cite{Bonn12}. Utilizing this method based on the fibration, we can deduce the heat kernel for $\Delta_{G_0}$ on $G_0$, denoted by $p_t^{G_0}$, as a subelliptic heat kernel on $G_0$. See Proposition~{\upshape\ref{TheHeatKernelOnASemiSimpleGroup:prop5.3.1}}. In the present subsection we demonstrate the connection between the two expressions $p_t^{G_0}$ and $\rho_t^{G_0}$ of the heat kernel on $G_0$.\\

It is clear that the Laplace--Beltrami operator $\Delta_{G_0}$ is a left $G_0$-invariant subelliptic second order differential operator. Hence we can deduce Proposition~{\upshape\ref{TheHeatKernelOnASemiSimpleGroup:prop5.3.1}} concerning an integral presentation of the heat kernel on $G_0$, which is established for a sub-Riemannian structure on $G_0$ in \cite{Bonn12}.

\begin{prop}[{Cf. \cite[Proposition 3.3]{Bonn12} for {$4\omega_{G_0}+\Delta_{K_0}$}}]\label{TheHeatKernelOnASemiSimpleGroup:prop5.3.1}

The heat kernel $p_t^{G_0}=p_t^{G_0}(r,\varphi)$ for $\Delta_{G_0}=\omega_{G_0}+\Delta_{K_0}$ on $G_0$ is given for $t>0$, $r>0$ and $\varphi\in\mathbb{R}$ by
\begin{equation}\label{TheHeatKernelOnASemiSimpleGroup:eq5.3.1}
\begin{split}
&p_t^{G_0}(r,\varphi)\\
&=\frac{2\pi e^{-\frac{t}{8}}}{(2\pi t)^2}\sum_{k\in\mathbb{Z}}\int_{\mathbb{R}}e^{-\frac{4\operatorname{Arc}\!\cosh^2\left(\cosh\frac{r}{2}\cosh\frac{y}{2}\right)-\left(\frac{y-i\varphi}{2}-2\pi ik\right)^2}{2t}}\frac{\operatorname{Arc}\!\cosh\left(\cosh\frac{r}{2}\cosh\frac{y}{2}\right)}{\sqrt{\left(\cosh\frac{r}{2}\cosh\frac{y}{2}\right)^2-1}}dy.
\end{split}
\end{equation}
\end{prop}
\noindent Here the difference of the factor $2\pi$ of \eqref{TheHeatKernelOnASemiSimpleGroup:eq5.3.1} arises from the fact that we take the normalized volume $\int_{K_0}dk=1$.

Let $g=k_{\theta_1/2}a_{r/2}k_{\theta_2/2}\in G_0$ be the Cartan decomposition with $\theta_1,\theta_2\in\mathbb{R}$ and $r\geq0$. Then we state the main result of this subsection.

\begin{prop}[The heat kernel on $G_0=SL(2,\mathbb{R})$]\label{TheHeatKernelOnASemiSimpleGroup:prop5.3.2}
For $\theta_1,\theta_2\in\mathbb{R}$ and $r\geq0$, it holds that
\begin{equation*}
\rho_t^{G_0}(k_{\theta_1/2}a_{r/2}k_{\theta_2/2})=p_t^{G_0}\left(\frac{r}{2},-\frac{\theta_1+\theta_2}{2}\right).
\end{equation*}
\end{prop}

We introduce the analytic continuation of $\mathbb{R}\ni\theta/2\mapsto k_{\theta/2}\in K_0$ given by
\begin{equation*}
k_{(\theta+iy)/2}:=
\begin{pmatrix}
\cos\frac{\theta+iy}{2} & -\sin\frac{\theta+iy}{2}\\
\sin\frac{\theta+iy}{2} & \cos\frac{\theta+iy}{2}
\end{pmatrix}\in G.
\end{equation*}
We then define the analytic continuation of $\rho_t^{K_0}$ as
\begin{equation*}
\tilde{\rho}_t^{K_0}(k_{(\theta+iy)/2}):=\frac{1}{4\pi}\vartheta\left(-\DS\frac{\theta+iy}{4\pi},\DS\frac{t}{2\pi}i\right).
\end{equation*}
Let $K^{\mathrm{ac}}_0:=\{k_z\ ;\ z\in\mathbb{C}\}$ be a subgroup of $G^{\tilde{\sigma}_0}$. Using the anti-generalized-Cartan-decomposition
\begin{equation*}
G=G^{\sigma}\overline{A^{+}}G^{\tilde{\sigma}_0},
\end{equation*}
we obtain a quotient space $G\times K^{\mathrm{ac}}_0/\operatorname{diag}K^{\mathrm{ac}}_0\approx G$.

\begin{proof}
By Fourier series expansion, we get an inversion formula for Jacobi's theta function $\vartheta$
\begin{equation*}
\frac{1}{\sqrt{2\pi t}}\sum_{k\in\mathbb{Z}}e^{-\frac{(z-2k\pi)^2}{2t}}=\frac{1}{2\pi}\sum_{k\in\mathbb{Z}}e^{-\frac{k^2t}{2}+ikz},
\end{equation*}
where $t>0$ and $z\in\mathbb{C}$. Hence we get
\begin{equation*}
\begin{split}
&p_t^{G_0}\left(\frac{r}{2},-\frac{\theta_1+\theta_2}{2}\right)\\
&=\frac{2\pi e^{-\frac{t}{8}}}{(2\pi t)^2}\int_{\mathbb{R}}e^{-\frac{4\operatorname{Arc}\!\cosh^2\left(\cosh\frac{r}{2}\cosh\frac{y}{2}\right)}{2t}}\frac{\operatorname{Arc}\!\cosh\left(\cosh\frac{r}{2}\cosh\frac{y}{2}\right)}{\sqrt{\left(\cosh\frac{r}{2}\cosh\frac{y}{2}\right)^2-1}}\sum_{k\in\mathbb{Z}}e^{-\frac{\left(\frac{(\theta_1+\theta_2)+iy}{2}+2k\pi\right)^2}{2t}}dy\\
&=\frac{e^{-\frac{t}{8}}}{(2\pi t)^{3/2}}\int_{\mathbb{R}}e^{-\frac{4\operatorname{Arc}\!\cosh^2\left(\cosh\frac{r}{2}\cosh\frac{y}{2}\right)}{2t}}\frac{\operatorname{Arc}\!\cosh\left(\cosh\frac{r}{2}\cosh\frac{y}{2}\right)}{\sqrt{\left(\cosh\frac{r}{2}\cosh\frac{y}{2}\right)^2-1}}\sum_{k\in\mathbb{Z}}e^{-\frac{k^2t}{2}- ik\frac{(\theta_1+\theta_2)+iy}{2}}dy.
\end{split}
\end{equation*}
Note that, for each $y\neq0$, the value $\operatorname{Arc}\!\cosh(\cosh(r/2)\cosh(y/2))$ can be regarded as the length of a right hyperbolic triangle with the legs of length $r/2$ and $|y|/2$ by the hyperbolic Pythagorean theorem. Changing the independent variables
\begin{equation*}
\cosh(s/2)=\cosh(r/2)\cosh(y/2),
\end{equation*}
we get
\begin{equation*}
\begin{split}
&\frac{e^{-\frac{t}{8}}}{(2\pi t)^{3/2}}\int_{\mathbb{R}}e^{-\frac{4\operatorname{Arc}\!\cosh^2\left(\cosh\frac{r}{2}\cosh\frac{y}{2}\right)}{2t}}\frac{\operatorname{Arc}\!\cosh\left(\cosh\frac{r}{2}\cosh\frac{y}{2}\right)}{\sqrt{\left(\cosh\frac{r}{2}\cosh\frac{y}{2}\right)^2-1}}dy\\
&=\frac{\sqrt{2}e^{-\frac{t}{8}}}{(2\pi t)^{3/2}}\int^{\infty}_r\frac{s e^{-\frac{s^2}{2t}}}{\sqrt{\cosh s-\cosh r}}ds=g_{t/2}^{G_0}(a_{r/2})=(\mathbf{M}g_{t/4}^G)(a_{r/2})\\
&=\frac{1}{2}\int^{\infty}_{s=\frac{r}{2}}\frac{g_{t/4}^G(a_{s/2})}{\sqrt{\cosh(2s)-\cosh r}}\frac{d(\cosh^2 s)}{\cosh s}=\frac{1}{4}\int^{\infty}_{s=r}\frac{g_{t/4}^G(a_{s/4})}{\sqrt{\cosh s-\cosh r}}\frac{d(\cosh^2 \frac{s}{2})}{\cosh \frac{s}{2}}\\
&=\frac{1}{8}\int^{\infty}_{y=-\infty}\frac{g_{t/4}^G(a_{2^{-1}\operatorname{Arc}\!\cosh(\cosh(r/2)\cosh(y/2))})}{\sqrt{\left(\cosh\frac{r}{2}\cosh\frac{y}{2}\right)^2-1}}\frac{d(\cosh^2(r/2)\cosh^2(y/2))}{\cosh(r/2)\cosh(y/2)}.
\end{split}
\end{equation*}
For $g\in G$ with $\operatorname{Iw}_A(g)=a_{r/2}$ and $0\neq y\in\mathbb{R}$, we set
\begin{equation*}
h_t^G(g,a_{y/2}):=\frac{1}{8}\frac{g_{t/4}^G(a_{2^{-1}\operatorname{Arc}\!\cosh(\cosh(r/2)\cosh(y/2))})}{\sqrt{\left(\cosh\frac{r}{2}\cosh\frac{y}{2}\right)^2-1}}.
\end{equation*}
Then the $K_0$-isotypic decomposition of $p_t^{G_0}$, cf. the proof of Lemma~{\upshape\ref{Chap2:lem2.1.1}}, yields
\begin{equation*}
\begin{split}
&p_t^{G_0}\left(\frac{r}{2},-\frac{\theta_1+\theta_2}{2}\right)\\
&=\sum_{k\in\mathbb{Z}}\int_{\mathbb{R}}\frac{e^{-\frac{t}{8}}}{(2\pi t)^{3/2}}e^{-\frac{4\operatorname{Arc}\!\cosh^2\left(\cosh\frac{r}{2}\cosh\frac{y}{2}\right)}{2t}}\frac{\operatorname{Arc}\!\cosh\left(\cosh\frac{r}{2}\cosh\frac{y}{2}\right)}{\sqrt{\left(\cosh\frac{r}{2}\cosh\frac{y}{2}\right)^2-1}}e^{-\frac{k^2t}{2}- ik\frac{(\theta_1+\theta_2)+iy}{2}}dy\\
&=\sum_{k\in\mathbb{Z}}\int^{\infty}_{y=-\infty}h_t^G(a_{r/2},a_{y/2})e^{-\frac{k^2t}{2}- ik\frac{(\theta_1+\theta_2)+iy}{2}}\frac{d(\cosh^2(r/2)\cosh^2(y/2))}{\cosh(r/2)\cosh(y/2)}\\
&=\int^{\infty}_{y=-\infty}h_t^G(a_{r/2},a_{y/2})\tilde{\rho}_t^{K_0}(k_{-\frac{(\theta_1+\theta_2)+iy}{2}})\frac{d(\cosh^2(r/2)\cosh^2(y/2))}{\cosh(r/2)\cosh(y/2)}
\end{split}
\end{equation*}
as functions on $G\times K^{\mathrm{ac}}_0/\operatorname{diag}K^{\mathrm{ac}}_0$. Since we have
\begin{equation*}
a_{r/2}k_{iy/2}=a_{r/2}\begin{pmatrix}
\cosh\frac{y}{2} & i\sinh\frac{y}{2}\\
-i\sinh\frac{y}{2} & \cosh\frac{y}{2}
\end{pmatrix}\in \exp\mathfrak{p},
\end{equation*}
it follows that, for each $r\geq0$ and $y\in\mathbb{R}$, there exists $x_{r,y}\in K$ such that
\begin{equation*}
a_{r/2}k_{iy/2}=x_{r,y}a_{\operatorname{Arc}\!\cosh(\cosh(r/2)\cosh(y/2))}x_{r,y}^{-1}\in\operatorname{Ad}_G(K)\overline{A^{+}}=\exp\mathfrak{p}.
\end{equation*}
Hence, as functions on $G\times K/\operatorname{diag}K$, it follows that
\begin{equation*}
\begin{split}
\int^{\infty}_{y=-\infty}&h_t^G(a_{r/2},a_{y/2})\tilde{\rho}_t^{K_0}(k_{-\frac{(\theta_1+\theta_2)+iy}{2}})\frac{d(\cosh^2(r/2)\cosh^2(y/2))}{\cosh(r/2)\cosh(y/2)}\\
=&\left[\int^{\infty}_{y=-\infty}h_t^G(\bullet,a_{y/2})\frac{d(\cosh^2(r/2)\cosh^2(y/2))}{\cosh(r/2)\cosh(y/2)}\times\rho_t^{K_0}\right]([a_{r/2}k_{iy/2},k_{-\frac{\theta_1+\theta_2}{2}}])\\
=&\left[\int^{\infty}_{y=-\infty}h_t^G(\bullet,a_{y/2})\frac{d(\cosh^2(r/2)\cosh^2(y/2))}{\cosh(r/2)\cosh(y/2)}\times\rho_t^{K_0}\right]\\
&([x_{r,y}a_{\operatorname{Arc}\!\cosh(\cosh(r/2)\cosh(y/2))}x_{r,y}^{-1},k_{-\frac{\theta_1+\theta_2}{2}}])\\
=&\left[\int^{\infty}_{y=-\infty}h_t^G(\bullet,e)\frac{d(\cosh^2(r/2)\cosh^2(y/2))}{\cosh(r/2)\cosh(y/2)}\times\rho_t^{K_0}\right]\\
&([a_{\operatorname{Arc}\!\cosh(\cosh(r/2)\cosh(y/2))},k_{-\frac{\theta_1+\theta_2}{2}}]).
\end{split}
\end{equation*}
Thus we deduce the identity of functions on $G_0\times K_0/\operatorname{diag}K_0$
\begin{equation*}
\begin{split}
&p_t^{G_0}\left(\frac{r}{2},-\frac{\theta_1+\theta_2}{2}\right)\\
&=\left[\int^{\infty}_{y=-\infty}h_t^G(\bullet,e)\frac{d(\cosh^2(r/2)\cosh^2(y/2))}{\cosh(r/2)\cosh(y/2)}\times\rho_t^{K_0}\right]([a_{\operatorname{Arc}\!\cosh(\cosh(r/2)\cosh(y/2))},k_{-\frac{\theta_1+\theta_2}{2}}])\\
&=\left[\frac{1}{2}\int^{\infty}_{s=\frac{r}{2}}\frac{g_{t/4}^G(\bullet)}{\sqrt{\cosh(2s)-\cosh r}}\frac{d(\cosh^2 s)}{\cosh s}\times\rho_t^{K_0}\right]([a_{s/2},k_{-\frac{\theta_1+\theta_2}{2}}])\\
&=[\mathbf{M}g_{t/4}^G\times\rho_t^{K_0}]([a_{r/2},k_{-\frac{\theta_1+\theta_2}{2}}])=[g_{t/2}^{G_0}\times\rho_t^{K_0}]([a_{r/2},k_{-\frac{\theta_1+\theta_2}{2}}])=\rho_t^{G_0}(k_{\theta_1/2}a_{r/2}k_{\theta_2/2}).
\end{split}
\end{equation*}
This concludes the proof.
\end{proof}




\bibliographystyle{abbrvurl}

\bibliography{TheHeatKernelOnASemisimpleLieGroupWithFiniteCenter}

\end{document}